\documentclass[12pt,a4paper]{amsart}
\usepackage[all]{xy}
\SelectTips{cm}{}
\calclayout

\usepackage{amscd,amssymb,amsopn,amsmath,amsthm,graphics,mathrsfs,accents}
\usepackage[colorlinks=true,linkcolor=red,citecolor=blue]{hyperref}
\usepackage[all]{xy}

\newcommand{\new}[1]{{\color{blue} #1}}

\newcommand{\rt}{\rightarrow}
\newcommand{\lrt}{\longrightarrow}

\def\xx{\text{\boldmath $x$}}

\def\yy{\text{\boldmath $y$}}
\newcommand{\N}{\mathbb{N} }

\def\MCM{\mathsf{MCM}}

\newcommand{\X} {\mathcal{X}}
\newcommand{\Md}{\mathsf{Mod}}
\newcommand{\md}{\mathsf{mod}}

\newcommand{\st}{\stackrel}

\newcommand{\im}{{\rm{Im}}}
\newcommand{\id}{{\rm{id}}}
\newcommand{\pd}{{\rm{pd}}}
\newcommand{\add} {\mathsf{add}}
\newcommand{\m}{{\mathfrak{m}}}

\newcommand{\p}{{\mathfrak{p}}}
\newcommand{\cok}{{\rm{Coker}}}

\newcommand{\ind}{{\rm{ind}}}

\newcommand{\HT}{\mathsf{H}}
\newcommand{\CM}{\mathsf{CM}}

\newcommand{\Hom}{{\mathsf{Hom}}}
\newcommand{\End}{{\mathsf{End}}}
\newcommand{\Add} {\mathsf{Add}}

\newcommand{\A} {\mathcal{A}}
\newcommand{\D} {\mathcal{D}}

\newcommand{\Ext}{\mathsf{{Ext}}}
\newcommand{\F}{\mathcal{{F}}}

\newcommand{\C}{\mathcal{{C}}}
\newcommand{\uh}{\underline{\mathsf{h}}}
\newcommand{\uHom}{\overline{\mathsf{Hom}}}
\newcommand{\uEnd}{\overline{\mathsf{End}}}
\newcommand{\Supp}{{\mathsf{supp}}}

\newcommand{\rad}{\rm{rad}}
\def\fd{\operatorname{\mathsf{FD}}}
\newtheorem{theorem}{Theorem}[section]
\newtheorem{cor}[theorem]{Corollary}

\newtheorem{lemma}[theorem]{Lemma}
\newtheorem{prop}[theorem]{Proposition}

\theoremstyle{definition}

\newtheorem{dfn}[theorem]{Definition}

\newtheorem{remark}[theorem]{Remark}
\newtheorem{s}[theorem]{}

\theoremstyle{plain}

\theoremstyle{definition}

\newtheorem{conv}[theorem]{Convention}

\numberwithin{equation}{section}

\begin{document}

\title[Representation-theoretic properties of balanced big $\CM$ modules] 
{Representation-theoretic properties of balanced big Cohen-Macaulay modules}

\author[Bahlekeh, Fotouhi and Salarian]{Abdolnaser Bahlekeh, Fahimeh Sadat Fotouhi and Shokrollah Salarian}

\address{Department of Mathematics, Gonbad Kavous University, Postal Code:49717-99151, Gonbad Kavous, Iran}
\email{bahlekeh@gonbad.ac.ir}

\address{Department of Mathematics, University of Isfahan, P.O.Box: 81746-73441, Isfahan, Iran}
\email{fsfotuhi@sci.ui.ac.ir}

\address{Department of Mathematics, University of Isfahan, P.O.Box: 81746-73441, Isfahan,
 Iran and \\ School of Mathematics, Institute for Research in Fundamental Science (IPM), P.O.Box: 19395-5746, Tehran, Iran}
 \email{Salarian@ipm.ir}

\subjclass[2000]{13C14, 16G60, 13H10, 13D07, 16E65.}

\keywords{balanced big Cohen-Macaulay modules, finite Cohen-Macaulay type, first Brauer-Thrall conjecture,
Gabriel-Roiter (co)measure, maximal Cohen-Macaulay modules, $\uh$-length.}

\thanks{The third author was partly supported by a grant from IPM (No. 95130218).\\
This paper is dedicated to the memory of the third author's brother, Rahmatollah Salarian,
who passed away while this paper was in preparation.}

\begin{abstract}
Let $(R, \m, k)$ be a complete Cohen-Macaulay local ring.
In this paper, we assign a numerical invariant, for any balanced big Cohen-Macaulay
module, called $\uh$-length. Among other results, it is proved that, for a given balanced big Cohen-Macaulay
$R$-module $M$ with an $\m$-primary cohomological annihilator, if there is a bound on the $\uh$-length
of all modules appearing in $\CM$-support of $M$,
then it is fully decomposable, i.e. it is a direct sum of finitely generated modules.
While the first Brauer-Thrall conjecture fails in general
by a counterexample of Dieterich dealing with multiplicities to measure the size of
maximal Cohen-Macaulay modules,
our formalism establishes the validity  of  the conjecture for complete Cohen-Macaulay local rings.
In addition, the pure-semisimplicity of  a subcategory of balanced big Cohen-Macaulay modules
is settled. Namely, it is shown that $R$ is of finite $\CM$-type if and only if
 $R$ is an isolated singularity and the category
of all fully decomposable  balanced big Cohen-Macaulay
modules  is  closed under kernels of epimorphisms.
Finally, we examine the mentioned results in the context of Cohen-Macaulay
artin algebras admitting a dualizing bimodule $\omega$, as defined by Auslander and Reiten.
It will turn out that, $\omega$-Gorenstein projective modules with bounded $\CM$-support
are fully decomposable.
In particular, a Cohen-Macaulay algebra $\Lambda$ is of finite $\CM$-type if and only if
every $\omega$-Gorenstein projective module is of finite $\CM$-type, which generalizes a result
of Chen for Gorenstein algebras.
Our main tool in the proof of results is Gabriel-Roiter
(co)measure, an invariant assigned to modules of finite length,
and defined by Gabriel and Ringel. This, in fact, provides an application of the  Gabriel-Roiter
(co)measure in the category of maximal Cohen-Macaulay modules.
\end{abstract}

\maketitle

\tableofcontents

\section{Introduction}

In representation theory of artin algebras, there is a large body of work
on the connections between representation-theoretic
properties of the category of finitely generated $\Lambda$-modules and global
structural properties of the algebra $\Lambda$. In this direction, the first Brauer-Thrall conjecture
asserts that if a finite-dimensional algebra $A$
over a field $k$ is of bounded representation type (meaning that
there is a bound on the length of the indecomposable finitely generated $A$-modules),
then $A$ is of finite representation type, i.e. the set of isomorphism classes
of indecomposable finitely generated modules is finite; see \cite{j}.
This  conjecture was proved by Roiter \cite{ro} and it is
proved by Ringel \cite{ri2, ri44} over artin algebras.
Another instance of this connection is the {\em pure-semisimple conjecture}
which predicts that every left pure-semisimple ring (a ring over which every left module
is a direct sum of finitely generated ones) is of finite representation type.
Left pure-semisimple rings are known to be left artinian by a result of Chase \cite[Theorem 4.4]{c}.
The validity of the pure-semisimple conjecture for artin algebras comes from a
famous result of Auslander \cite{a1, a11} (see also Ringel-Tachikawa \cite[Corollary 4.4]{rt}),
where they have shown that an artin algebra $\Lambda$ is of finite representation
type if and only if every left $\Lambda$-module is a direct sum of finitely generated
modules. Motivated by Auslander's result, studying decomposition of Gorenstein projective
modules over artin algebras into finitely generated ones has been the subject of several expositions
(see \cite{be, che, lz, ru, ru1}). In particular, a result of Beligiannis \cite[Theorem 4.10]{be}
asserts that a virtually Gorenstein algebra $\Lambda$ is of finite Gorenstein representation type,
in the sense that there are only finitely many isomorphism classes of indecomposable
finitely generated Gorenstein projective $\Lambda$-modules if and only if
 any left Gorenstein projective
$\Lambda$-module is a direct sum of finitely generated ones.
This solves a problem raised by Chen \cite{che}, who proved it for
Gorenstein artin algebras.

On the other hand, over the past several decades Cohen-Macaulay rings
and maximal Cohen-Macaulay modules have achieved a great deal of significance
in commutative algebra and algebraic geometry.
Hochster and Huneke \cite{hh} write that for many theorems ``the Cohen-Macaulay condition
(possibly on the local rings of a variety) is just what is needed
to make the theory work.''
Let $(R, \m, k)$ be a commutative noetherian local ring. Hochster \cite{hoc}
defines a not necessarily finitely generated $R$-module
$M$ is {\em big Cohen-Macaulay}, if there exists a system of parameters of $R$
which is an $M$-regular sequence. Sharp \cite{sha} called a big Cohen-Macaulay
$R$-module $M$ is {\em (weak) balanced big Cohen-Macaulay,}
((weak) balanced big $\CM$, for short), provided that every system of parameters of $R$
is an (a weak) $M$-regular sequence. A finitely generated $R$-module $M$
is maximal Cohen-Macaulay (abbreviated, $\MCM$), if it is either balanced big Cohen-Macaulay
or zero.

Motivated by the above mentioned results, the major issues considered in this
paper are when a given balanced big $\CM$ module is a direct
sum of finitely generated modules; when every balanced big $\CM$ module
is so; analogues of the first Brauer-Thrall conjecture for modules
and analogues result for $\omega$-Gorenstein projective modules over Cohen-Macaulay
artin algebras in the sense of Auslander and Reiten \cite{AR1, AR}.

A natural interpretation of the first Brauer-Thrall conjecture
in this context, states that a commutative noetherian local ring $(R, \m)$
is of finite Cohen-Macaulay type, provided that there is a bound on the multiplicities
of indecomposable $\MCM$ modules. Recall that $R$ is said to be of
finite Cohen-Macaulay type (finite $\CM$-type, for short), if there are only
finitely many non-isomorphic indecomposable $\MCM$ $R$-modules.
An example discovered by Dieterich \cite{di}, disproved the conjecture in general.
However, over several classes of rings, this conjecture is known to be true.
Namely, it has been answered affirmatively for complete, equicharacteristic
Cohen-Macaulay isolated singularities over a perfect field, independently by
Dietrich \cite{di} and Yoshino \cite{yo1}. This result was extended by Leuschke and
Wiegand \cite[Theorem 3.4]{lw2} to the case where the ring is equicharacteristic excellent
with algebraically closed residue field $k$.
On the other hand, inspired by the pure-semisimplicity conjecture,
Beligiannis \cite[Theorem 4.20]{be} has shown that
a commutative noetherian Gorenstein complete local ring $R$
being of finite $\CM$-type is tantamount to saying that
any Gorenstein projective $R$-module is a direct sum of finitely generated modules.

In this paper, we focus our attention on modules of finite type.
In fact, we will treat the support of a module,
instead of all finitely generated indecomposable modules.
Recall that the {\em support} of a module $M$ over an artin algebra $\Lambda$,
denoted by $\Supp_{\Lambda}(M)$, is the set of all indecomposable finitely generated $\Lambda$-modules
$N$ such that $\Hom_{\Lambda}(N,M)\neq 0$. It is a consequence of nice
results of Auslander \cite[Theorem B]{a2} and also Ringel \cite[Theorem 1]{ri44} that, for a given $\Lambda$-module
$M$, if $\Supp_{\Lambda}(M)$ is of bounded representation type (meaning that
there is a bound on the length of modules in $\Supp_{\Lambda}(M)$),
then $M$ is of finite type. Recall that a $\Lambda$-module $M$
is said to be of finite type, provided
it is the direct sum of (arbitrarily many) copies of a finite
number, up to isomorphism, of indecomposable modules of finite length; see \cite{ri44}.
The main tool in Ringel's proof is Gabriel-Roiter
(co)measure, an invariant assigned to any module of finite length,
and defined by Gabriel and Ringel \cite{ga, ri3, ri1} based on
Roiter's induction scheme in his proof of the first Brauer-Thrall conjecture.

In order to state our results precisely, let us recall some notions.

From now on, assume that $(R, \m, k)$ is a commutative noetherian complete
Cohen-Macaulay local ring with a canonical module $\omega$.
We say that a balanced big $\CM$ $R$-module $M$ is {\em of finite
$\CM$-type}, if it is a direct sum of (arbitrarily many)
copies of a finite number, up to isomorphisms, of indecomposable $\MCM$ modules
and it is said to be {\em fully decomposable},
provided it is a direct sum of finitely generated modules.
The class of all fully decomposable modules will be denoted by $\fd$.

Moreover, by $\CM$-support of a balanced big $\CM$ $R$-module $M$, denoted by $\CM$-$\Supp_{R}(M)$,
we mean the set of all indecomposable $\MCM$ $R$-modules
$N$ such that $\Hom_R(N,M)\neq 0$.
For a (not necessarily finitely generated) balanced big $\CM$ $R$-module $M$,
we set $\uh(M)=\underline{\Hom}_R(M,M\oplus G)$, where  $\alpha:G\lrt k$ is a right minimal $\MCM$-approximation.
We say that $M$ has finite {\em $\uh$-length,} provided that $l_R(\uh(M))<\infty$.
Also, $M$ is said to have an {\em $\m$-primary cohomological annihilator}, if $\m^t\uh(M)=0$,
for $t\gg 0$. One should observe that, this is equivalent to saying that
$\m^t\Ext_R^1(M, -)=0$, by a theorem of Hilton-Rees \cite{hr}.

Section 2 of the paper, is devoted  to comparing the length of the stable $\Hom$ and $\uh$-length
of maximal Cohen-Macaulay modules with classical invariants such as multiplicity and Betti number.

The main result in section 3 enables us to demonstrate the utilization of the
Gabriel-Roiter (co)measure for the category of balanced big Cohen-Macaulay modules; see Theorem \ref{the2}.
The purpose of section 4 is to study balanced big $\CM$ modules with bounded $\CM$-support.
In particular, we prove the result below; see Theorems \ref{th7} and \ref{th2}.

\begin{theorem}A balanced big $\CM$ $R$-module having an $\m$-primary
cohomological annihilator with bounded $\uh$-length on $\CM$-support is fully decomposable.
In particular, any balanced big $\CM$ modules with an $\m$-primary
cohomological annihilator and of bounded $\uh$-length on $\CM$-support, satisfies complements direct summands.
\end{theorem}

In section 5, we investigate balanced big $\CM$ modules  with  large
(finite) $\uh$-length, for instance, we have  the following  result; see Theorem \ref{th1}.

\begin{theorem}\label{ddim}Let $R$ be an isolated singularity and
let $M$ be a balanced big $\CM$ $R$-module with an $\m$-primary cohomological annihilator.
If $M$ is not of finite $\CM$-type, then 
there are indecomposable $\MCM$ $R$-modules of arbitrarily large
(finite) $\uh$-length.
\end{theorem}
It should be noted that this result provides a kind of the first Brauer-Thrall theorem for modules.
In particular, it guarantees the validity of the first Brauer-Thrall
conjecture for complete Cohen-Macaulay local rings, considering $\uh$-length
as an invariant to measure the size of $\MCM$ modules. Indeed, we have the
result below; see Corollary \ref{cor2}.

\begin{cor}Let the category of all indecomposable $\MCM$ $R$-modules
be of bounded $\uh$-length. Then $R$ is of finite $\CM$-type.
\end{cor}
We would like to point out that, as already mentioned previously,
the first Brauer-Thrall conjecture fails in general when multiplicity
is used as the size, by an example of Dieterich \cite{di}.

In addition, it will be observed that the representation-theoretic properties of balanced big $\CM$ modules
have important consequences for the structural properties
of the ring. Actually, Theorem \ref{prop7} asserts that:
\begin{theorem}If any balanced big $\CM$
$R$-module $M$ admitting a right resolution by modules in $\Add\omega$,
is fully decomposable, then $R$ is an isolated singularity.
\end{theorem}
It seems that this result is a generalization of a result of Chase \cite{c}
for the category of $\MCM$ modules. Furthermore,
we prove a variant of a celebrated theorem of Auslander \cite{a1, a11},
Ringel-Tachikawa \cite{rt}, Chen \cite{che} and Beligiannis \cite{be} for Cohen-Macaulay local rings.
In fact, our main result in section 6 reads as follows.
\begin{theorem}\label{fcmti}A complete Cohen-Macaulay local ring $R$
is of finite $\CM$-type if and only if the category of balanced big $\CM$ $R$-modules
with $\m$-primary cohomological annihilators coincides with the category of
fully decomposable balanced big $\CM$ modules. Equivalently; $R$ is an isolated singularity and
the category of all fully decomposable balanced big $\CM$
modules is closed under kernels of epimorphisms.
\end{theorem}
The precise statement of the above result is Theorem \ref{th6}.

In the paper's final section, we are concerned with Cohen-Macaulay
artin algebras and Cohen-Macaulay modules in the sense of Auslander and Reiten \cite{AR, AR1}.
Recall that an artin algebra $\Lambda$ is said to be a {\em Cohen-Macaulay algebra},
if there is a pair of adjoint functors $(G, F)$ on the category of finitely
generated (left) $\Lambda$-module, $\md\Lambda$,
which induce mutually inverse equivalences between the full subcategories of
$\md\Lambda$ consisting of the $\Lambda$-modules of finite injective dimension
and the $\Lambda$-modules of finite projective dimension.
It is known that an artin algebra $\Lambda$ is Cohen-Macaulay if and only if
there is a $\Lambda$-bimodule $\omega$ such that the pair of adjoint functors
$(\omega\otimes_{\Lambda}-, \Hom_{\Lambda}(\omega, -))$ has the desired properties.
In this case, $\omega$ is called a dualizing module for $\Lambda$.
A (not necessarily finitely generated) $\Lambda$-module $M$ is said to be
$\omega$-{\em{Gorenstein projective}}, provided that it admits a right resolution by modules in $\Add\omega$.
Following Auslander and Reiten \cite{AR1}, a finitely generated $\omega$-Gorenstein
projective module will be called a {\em Cohen-Macaulay} module.
The notion of Cohen-Macaulay artin algebras (and also Cohen-Macaulay modules) is generalizations of
commutative complete Cohen-Macaulay local rings as well as Gorenstein artin
algebras (Gorenstein projective modules). {Recall that an artin algebra
$\Lambda$ is said to be a {\em Gorenstein algebra}, provided the injective dimension of $_{\Lambda}\Lambda$
as well as of $\Lambda_{\Lambda}$ is finite.}
The main goal of section 7 is to study the decomposition
properties of $\omega$-Gorenstein projective modules in connection with the property
that $\Lambda$ is of finite $\CM$-type. In this direction, it is proved that any
$\omega$-Gorenstein projective $\Lambda$-module $M$ in which $\CM$-$\Supp_{\Lambda}(M)$
is of bounded length, is fully decomposable; see Theorem \ref{the9}.
Using this result, we prove that there exist indecomposable $\CM$ $\Lambda$-modules
of (arbitrarily) large finite length, if there is an $\omega$-Gorenstein projective
$\Lambda$-module which is not of finite $\CM$-type; see Theorem \ref{th4}. This is fruitful from the
point of view that it is an analog of the first Brauer-Thrall theorem for modules over Cohen-Macaulay artin
algebras; see Corollary \ref{cor10}.
In addition, we extend Chen's result \cite[Main theorem]{che}
to Cohen-Macaulay artin algebras.
Namely, it is shown that $\Lambda$ is of finite $\CM$-type if and only if
every $\omega$-Gorenstein projective module is fully decomposable;
see Theorem \ref{th3}.

We would like to emphasize that in proving our results,
we strongly use the notion of Gabriel-Roiter (co)measure;
see \ref{s1} for the definition of Gabriel-Roiter (co)measure.
So our method is totally different from the previous ones
which are based on functorial approach; see for example \cite{{be, che}}.
It is well understood that the Gabriel-Roiter (co)measure is a helpful invariant
dealing with representations of an artin algebra; see \cite{ri3, ri1, ch}.
So it seems worthwhile to unfold the use of this notion in the setting of
commutative noetherian rings. In this direction, our point of view gives
another nice feature of the paper which brings the use
of Gabriel-Roiter (co)measure in the context of $\MCM$
modules; see also Theorem \ref{the2}.

\section{Preliminary results}
This section is devoted to stating the definitions and basic properties
of notions which we will freely use in the later sections.
We also define length of the stable $\Hom$, $\uh$-length, of balanced big Cohen-Macaulay
modules and study its relationship with well-known invariants, such as multiplicity and Betti number.
Let us start with our convention.

\begin{conv}
Throughout the paper, unless otherwise specified, $(R, \m, k)$ is a $d$-dimensional commutative
complete Cohen-Macaulay local ring with a dualizing (or canonical) module $\omega$.
The category of all (finitely generated) $R$-modules will be denoted by
($\md R$) $\Md R$.
\end{conv}

\begin{s}An $R$-homomorphism $f:X\lrt Y$ is called {\em right minimal}, provided that any $R$-homomorphism
$g:X\lrt X$ satisfying $fg=f$, is an isomorphism.\\
 An $R$-homomorphism $f:M\lrt X$ with $M$ is $\MCM$ is called
a {\em right $\MCM$-approximation}, if the map $\Hom_R(L, f):\Hom_{R}(L, M)\lrt\Hom_R(L, X)$
is surjective for any $\MCM$ $R$-module $L$; and a right minimal
$\MCM$-approximation if, in addition, $f$ is right minimal.\\
It should be noted that by \cite[Theorem A]{ab}, every
finitely generated $R$-module admits a right minimal $\MCM$-approximation.
In the rest of this paper, we assume that $\alpha:G\lrt k$ is a right minimal
$\MCM$-approximation of the residue field $k$.
\end{s}

\begin{dfn}\label{bbcm}
(1) A local ring $(R,\m)$ is called an {\em isolated singularity},
if $R_\p$ is a regular ring for all nonmaximal prime ideals $\p$ of $R$.\\
(2) A finitely generated $R$-module $M$ is said to be {\em locally free
on the punctured spectrum of $R$}, if $M_\p$ is a free $R_\p$-module for all
nonmaximal prime ideals $\p$ of $R$.\\
(3) A system of parameters  $\xx=x_1, \cdots, x_d$ of $R$
is said to be a {\em faithful system of parameters}, if it annihilates $\Ext_R^1(M, N)$
for any $M$ in $\MCM$ modules and $N\in\md R$, (\cite[Definition 14.8]{lw1}).
If $R$-modules $M$ are taken from a subcategory $\C$ of $\MCM$ modules,
then we will say that $\xx$ is a faithful system of parameters for $\C$.
In the remainder of this paper, $\xx^t$, where $t>0$ is an integer,
stands for the ideal $(x_1^t, \cdots, x_d^t)$.\\
(4) Let $M$ be an $R$-module. A sequence of elements $\xx=x_1, \cdots, x_n\in\m$
is called a {\em weak $M$-regular sequence}, provided that $x_i$
is a non-zerodivisor on $M/{(x_1,\cdots, x_{i-1})M}$ for any $1\leq i\leq n$
(for $i=1$, we mean that $x_1$ is a non-zerodivisor on $M$). If, in addition,
$(x_1, \cdots, x_n)M\neq M$, then $\xx$ is said to be an {\em $M$-regular sequence.}
It is worth remarking that if $M$ is a non-zero finitely generated
$R$-module, then it follows from Nakayama's lemma that any weak $M$-regular
sequence is an $M$-regular sequence, as well. It is also known that
over local rings, any permutation of $M$-regular sequence, is again an $M$-regular sequence.\\
\end{dfn}

\begin{s}\label{s4} (1) We use $\X_{\omega}$ to denote the subcategory
consisting of all $R$-modules $M$ admitting
a right resolution by modules in $\Add\omega$, that is, an exact sequence of $R$-modules;
$$0\lrt M\lrt w_0\st{d_0}\lrt w_1\st{d_1}\lrt\cdots \st{d_{i-1}}\lrt w_i\st{d_i}\lrt\cdots,$$ with $w_i\in{\Add\omega}$.
By ${\Add}\omega$  (resp. $\add\omega$) we mean the full subcategory
of ${\Md}R$ (resp. $\md R$) consisting of all modules isomorphic to direct
summands of direct sums (resp. finite direct sums) of copies of $\omega$.
It is known that $\MCM=\X_{\omega}\bigcap\md R$. {To see this, according to \cite[Theorem 3.3.10]{bh},
a given module $M$ is $\MCM$ if and only if $\Ext_R^i(M, \omega)=0=\Ext_R^i(M^{*}, \omega)$ for all $i\geq 1$ and 
the natural homomorphism $\delta:M\lrt M^{**}$ is an isomorphism, where $M^{*}=\Hom_R(M, \omega)$.
Now assume that $M$ is an arbitrary $\MCM$ module and $\cdots\lrt P_1\lrt P_0\lrt M^{*}\lrt 0$
is an exact sequence in $\md R$ such that each $P_i$ is projective. So,  applying the functor $\Hom_R(-, \omega)$,
implies  that $M$ admits a right resolution by modules in $\add\omega$, giving the containment $\MCM\subseteq\X_{\omega}\bigcap\md R$.
For the opposite containment, take a finitely generated $R$-mdule $M$ in $\X_{\omega}$. Since $\omega$ has finite injective dimension and $\Ext_R^i(\omega, \omega)=0$ for all $i\geq 1$, one may deduce that $\Ext_R^i(M, \omega)=0$ for all $i\geq 1$. This, in turn, implies that $M^*\in\X_{\omega}$ and so  $\Ext_R^i(M^*, \omega)=0$ for all $i\geq 1$. 
Now, one may use the fact that $\delta_{\omega}:\omega\lrt\omega^{**}$ is an isomorphism, and conclude that the same is true for $\delta_M:M\lrt M^{**}$.
Hence $M$ will be a $\MCM$ module.} \\
(2) Recall that for a subcategory $\X$ of ${\md}R$, we let $\widehat{\X}$ denote the category whose objects are
the modules $M$ for which there is an exact sequence of $R$-modules;
$0\lrt X_n\lrt\cdots \lrt X_0\lrt M\lrt 0$ with $X_i\in\X$.
\end{s}

Now we introduce the notion of $\uh$-length, an invariant to measure the size of
balanced big $\CM$ $R$-modules.
\begin{dfn} $(i)$ For a given balanced big $\CM$ $R$-module $M$, we set $\uh(M)=\underline{\Hom}_R(M,M\oplus G)$
and  define {\em $\uh$-length} of $M$ as $l_R(\uh(M))$.

Assume that $R$ is an isolated singularity. So in view of \cite[Lemma 3.3]{yo},
${\underline{\Hom}}_R(M, N)$ is an artinian $R$-module, for all $\MCM$ $R$-modules $M$ and $N$.
In particular, any $\MCM$ $R$-module $M$ has finite $\uh$-length.
Recall that for any two $R$-modules $M$ and $N$,
$\underline{\Hom}_{R}(M, N)=\Hom_{R}(M, N)/{\mathfrak{P}(M, N)},$
where $\mathfrak{P}(M, N)$ is the $R$-submodule of $\Hom_{R}(M, N)$ consisting of all
homomorphisms factoring through projective modules.\\
$(ii)$ Let $\C$ be a subcategory of $\MCM$ $R$-modules.
We say that $\C$ is {\em of bounded $\uh$-length,}
if there is an integer $b>0$ such that
$|\uh(\C)|= {\text {sup}}\{l_R(\uh(M)) |M\in\C\}<b$.
\end{dfn}



\begin{prop}\label{cor1} Let $\C$ be a subcategory of $\MCM$ $R$-modules of
bounded $\uh$-length. Then there is a system of parameters $\xx$
 such that $\xx\uh(M)=0$ for all $M\in\C$. In particular, $\C$ admits a faithful system of parameters.
\end{prop}
\begin{proof}Take an integer $b>0$ such that $|\uh(\C)|<b$. So for any
$M\in\C$, $l_R(\uh(M))<b$ implying that $\m^b\uh(M)=0$.
Now choosing a system of parameters $\xx\in\m^b$, one gets that  $\xx{\uh}(M)=0$.
In particular, $\xx\underline{\Hom}_R(M,M)=0$, for any $M\in\C$.
Hence by a theorem of Hilton-Rees \cite{hr}, we infer that
$\xx\Ext_R^1(M, -)=0$. So the proof is finished.
\end{proof}

For a given finitely generated $R$-module $M$, by the {\em Betti number of
$M$,} $\beta(M)$, we mean the minimal number of generators for $M$.

\begin{lemma}\label{lem7}
Let $\C$ be a subcategory of $\MCM$ $R$-modules. Then
$\C$ has a bound on multiplicities if and only if it has a bound on
Betti numbers.
\end{lemma}
\begin{proof}
Assume that there exists an integer $b>0$ such that for any $M\in\C$, $\beta(M)<b$.
So for each $M\in\C$, there is an $R$-epimorphism $f:R^{b}\lrt M$.
Take a system of parameters $\xx=x_1, x_2, \cdots, x_d$ of $R$.
Tensoring $f$ with $R/{\xx R}$ over $R$,
gives rise to the epimorphism $\bar{f}:R^{b}/{\xx R^{b}}\lrt M/{\xx M}$, implying that
$l_R(M/{\xx M})<l_R(R/{\xx R})b$. According to
\cite[Proposition 1.7]{yo}, the multiplicity of $M$, $e(M)$, is less than or equal to
$l_R(M/{\xx M})$. Consequently, for any $M\in\C$, $e(M)<l_R(R/{\xx R})b$.
The other direction follows from the well-known fact that, for any $\MCM$ module
$M$, $\beta(M)$ is less than or equal to $e(M)$.
The proof then is completed.
\end{proof}

The next results show that there is a tight connection between the
invariants $\uh$-length and multiplicity of $\MCM$ modules.

\begin{lemma}\label{lem88}
Let $\C$ be a subcategory of $\MCM$ $R$-modules and let $\xx$ be a faithful system
of parameters for $\C$. If there is a bound on the multiplicities of modules in $\C$,
then $\C$ is of bounded $\uh$-length.
\end{lemma}
\begin{proof}
We first claim that there is an integer $b>0$ such that for any $M\in\C$, $l_R(\Hom_{R/\xx R}(M/\xx M,(M\oplus G)/\xx (M\oplus G)))<b$.
To do this, one should note that according to Lemma \ref{lem7}, there is a bound on the Betti numbers of modules in $\C$,
say $n$. Assume that $M$ is an arbitrary object of $\C$.
So there exist $R$-epimorphisms,
$f:R^n\lrt M$ and $g:R^n\lrt M\oplus G$. Tensoring $f$ and $g$ with
$R/\xx R$ over $R$, gives rise to the epimorphisms $\bar{f}:R^n/\xx R^n\lrt M/\xx M$ and
$\bar{g}:R^n/\xx R^n\lrt (M\oplus G)/\xx(M\oplus G)$. Now, the $R/{\xx R}$ (and also $R$)-monomorphism;
$$\Hom_{R/{\xx R}}(M/{\xx M}, (M\oplus G)/\xx(M\oplus G)) \lrt\Hom_{R/{\xx R}}(R^n/{\xx R^n}, (M\oplus G)/\xx(M\oplus G)),$$
together with  $R/{\xx R}$ (and also $R$)-epimorphism;
$$\Hom_{R/{\xx R}}(R^n/{\xx R^n}, R^n/{\xx R^n}) \lrt\Hom_{R/{\xx R}}(R^n/{\xx R^n}, (M\oplus G)/\xx(M\oplus G)),$$
lead us to obtain the inequality
$$l_R(\Hom_{R/\xx R}(M/{\xx M}, (M\oplus G)/\xx(M\oplus G)))\leq l_R(\Hom_{R/{\xx R}}(R^n/{\xx R^n}, R^n/{\xx R^n})),$$
giving the claim, because the right hand side is finite.
On the other hand, for any $0\leq i\leq d-1$, we have the following
exact sequence of $R$-modules; $$0\rt (M\oplus G)/{\xx_i(M\oplus G)}\st{x_{i+1}}\rt (M\oplus G)/{\xx_i(M\oplus G)}
\rt (M\oplus G)/\xx_{i+1}(M\oplus G)\rt 0,$$  where $\xx_i=x_1, \cdots,x_i$
(in case, $i=0$, we mean $\xx_i=0$).
This induces the exact sequence of $R$-modules; $\underline{\Hom}_R(M, (M\oplus G)/{\xx_i(M\oplus G)})\st{x_{i+1}}\lrt
\underline{\Hom}_R(M, (M\oplus G)/{\xx_i(M\oplus G)})\st{\phi}
\lrt\underline{\Hom}_R(M,M\oplus G/\xx_{i+1}(M\oplus G)).$
Since $\xx=\xx_d$ is a faithful system of parameters for $\C$ and $\underline{\Hom}_R(M,(M\oplus G)/{\xx_i(M\oplus G)})$ is a submodule of
$\Ext_R^1(M,\Omega^1_R((M\oplus G)/{\xx_i(M\oplus G)}))$, $\phi$ will be a monomorphism.
Now, it is easily seen that the inequality $l_R(\underline{\Hom}_R(M,(M\oplus G)))\leq l_R(\Hom_R(M, (M\oplus G)/\xx(M\oplus G)))$
holds true. By \cite[Lemma 2(ii) page 140]{ma},
we have an isomorphism $\Hom_R(M, (M\oplus G)/\xx(M\oplus G))\cong \Hom_{R/\xx R}(M/\xx M, (M\oplus G)/\xx(M\oplus G))$. Hence,
there is an integer $b>0$ such that $l_R(\underline{\Hom}_R(M, M\oplus G))<b$, as needed.
\end{proof}

\begin{lemma}\label{lem9}Let $\C$ be a subcategory of indecomposable
$\MCM$ $R$-modules. If $\C$ is of bounded $\uh$-length, then it has a bound
on multiplicies.
\end{lemma}
\begin{proof}
Take an integer $t>0$ such that for any object $M$ in $\C$, $l_R(\uh(M))<t$.
 Assume that $M$ is an arbitrary non-projective object of $\C$. Applying the functor
$\Hom_R(M, -)$ to the short exact sequence of $R$-modules; $0\lrt \m\lrt R\lrt k\lrt 0$,
gives rise to the following  exact sequence; $$0\lrt\Hom_R(M, \m)\lrt\Hom_R(M, R)\lrt\Hom_R(M, k)\lrt\underline{\Hom}_R(M, k)\lrt 0.$$
Since $M$ is non-projective, it can be easily seen that
the isomorphism of $R$-modules $\underline{\Hom}_R(M,k)\cong \Hom_R(M,k)$ holds true.
On the other hand, $\alpha:G\lrt k$ being a right minimal $\MCM$-approximation
forces $\underline{\Hom}_R(M,G)\lrt\underline{\Hom}_R(M,k)$ to be an epimorphism.
All of these facts together enable us to infer that $l_R(\Hom_R(M, k))<t$, because $l_R(\uh(M))<t$.
This indeed means that there is a bound on the Betti numbers of modules in $\C$.
Now Lemma \ref{lem7} gives the desired result.
\end{proof}
{\begin{remark}\label{re1} As we have mentioned in the introduction, a given $R$-module $M$ has an  $\m$-primary cohomological annihilator if and only if $\m^t\Ext_R^1(M, -)=0$ for some integer $t\gg 0$. Assume that $0\lrt M'\lrt M\lrt M''\lrt 0$ is a short exact sequence of $R$-modules. So by applying the functor $\Ext_R(-, N)$, where $N$ is an arbitrary $R$-module, one may infer that  the class consisting of modules with $\m$-primary cohomological annihilators, is closed under extensions
and kernels of epimorphisms.
\end{remark}
\begin{lemma}\label{re2} Let $f:M\lrt\oplus_{i\in I}M_i$ be an $R$-homomorphism, where each $M_i$ is finitely generated. If for {a sequence} $\xx=x_1, \cdots, x_n\in\m$, $\bar{f}:M/{\xx M}\lrt\oplus_{i\in I}M_i/{\xx M_i}$ is an epimorphism, then $f$ is so.
\end{lemma}
\begin{proof}
In order to obtain the desired result, it suffices to show that for any finite subset $J$ of $I$, the composition map $M\st{f}\lrt \oplus_{i\in I}M_i \st{h}\lrt\oplus_{i\in J}M_{i}$ is an epimorphism,
where $h$ is the projection map. Consider the composition map
$M/{\xx M}\st{f\otimes R/{\xx R}}\lrt \oplus_{i\in I}M_i/{\xx M_i}\st{h\otimes R/{\xx R}} \lrt\oplus_{i\in J}{M_i}/{\xx M_i}$, which evidently is an epimorphism. Thus letting $\cok{(h f)}=Z$, we have $Z/{\xx Z}=0$ and so Nakayama's lemma implies that $Z=0$, meaning that $hf$ is an epimorphism. Consequently, $f$ is an epimorphism, as needed.
\end{proof}

\begin{lemma}\label{re3}Let $0\lrt M'\lrt M\lrt M''\lrt 0$ be a pure exact sequence of $R$-modules. If $M', M$ are weak balanced big $\CM$ module, then $M''$ is so.
\end{lemma}
\begin{proof}
Assume that $\xx=x_1, \cdots,x_d$ is a system of parameters of $R$.
We must show that for any $1\leq i\leq d$, $x_i$ is a non-zerodivisor on $M''/{(x_1,\cdots,x_{i-1})M''}$. 
Since the sequence $0\lrt M'\lrt M\lrt M''\lrt 0$ is pure exact, tensoring this sequence with $R/{xR}$,
for any $x\in R$, gives us again a pure exact sequence. This fact allows us to prove the case only for $i=1$.
So assume that $x_1=x$ is a regular element of $R$.
Consider the following commutative diagram with exact rows; 
{\tiny{\footnotesize{$$\begin{CD}
0 @>>>  M'  @>>>  M @>>>  M'' @>>> 0\\
& & @V x VV @V x VV @V x VV & &\\ 0 @>>>  M' @>>>  M @>>>  M''
@>>> 0\\
& & @V VV @V VV @V VV & &\\ 0 @>>>  M'/{xM'} @>>>  M/{xM} @>>>  M''/{xM''}
@>>> 0.\end{CD}$$}}}  Since $M', M$ are weak balanced big $\CM$, the left multiplicative maps are
monomorphism. So, one may apply the snake lemma and deduce that the right multiplicative map
is also monomorphism, giving the desired result.
\end{proof}}

\section{Using Gabriel-Roiter (co)measure in the category of $\MCM$ modules}
This section is devoted to bring the use of Gabriel-Roiter (co)measure
in the category of $\MCM$ $R$-modules.
The notion of Gabriel-Roiter (co)measure, an invariant assigned to any
module of finite length, was defined by Gabriel and Ringel \cite{ga, ri3, ri1}.
Since this notion is a basic tool in proving the results of the paper,
we recall it and some of its properties which will be used later.

\begin{s}\label{s2}{\sc Gabriel-Roiter (co)measure:} Let $\Lambda$ be an artin algebra and
$M$ a finitely generated $\Lambda$-module.
The Gabriel-Roiter measure of $M$, denoted by $\mu(M)$,
was defined in \cite{ri3} by induction on the length of modules as follows: let
$\mu(0)=0$. Given a non-zero module $M$, we may assume by
induction that $\mu(N)$ is already defined for any proper
submodule $N$ of $M$. Set

\[ \mu (M)= {\rm max}\{\mu(N)\} + \left\{ {\begin{array}{ll}
0, & \text{ if $M$ is decomposable,}\\
\frac{1}{2^{l_{\Lambda}(M)}}, & \text{ if $M$ is indecomposable,}
\end{array}} \right. \]
here maximum is taken over all proper submodules $N$ of $M$ and $l_{\Lambda}(M)$
denotes the length of $M$ over $\Lambda$. Note that the maximum
always exists. We should refer the reader
to \cite{ri1} for an equivalent definition using subsets of natural numbers,
which reformulates Gabriel's definition. The Gabriel-Roiter comeasure of $M$,
denoted by $\mu^*(M)$, is defined as $-\mu(D(M))$, where $D(M)=\Hom_{\Lambda}(M, \coprod(E(S)))$
in which $E(S)$ runs over all injective envelope of
 simple $\Lambda$-modules.
\end{s}

\begin{s}\label{s1}
Let us make a list of several basic properties of Gabriel-Roiter (co)measure,
which have been proved by Ringel in \cite{ri3} and \cite{ri1}.\\
 {\sc Property 1}. {\em Let $Y$ be a $\Lambda$-module of finite length and $X\subseteq Y$ a submodule.
Then $\mu(X)\leq\mu(Y)$. If $Y$ is indecomposable and $X$ is a proper submodule $Y$, then $\mu(X)<\mu(Y)$.}\\
 {\sc Property 2}. {\em If $M$ is an indecomposable $\Lambda$-module of length $n$,
then $\mu(M)=a/{2^n}$ where $a$ is an odd natural number such that $2^{n-1}\leq a<2^n$}.\\
It follows from this property that any subcategory of $\ind(\md\Lambda)$ with bounded length
has only finitely many Gabriel-Roiter measures. Moreover, in view of
the equality $l_{\Lambda}(M)=l_{\Lambda}(D(M))$ and the definition of
Gabriel-Roiter comeasure, such a subcategory has also only
finitely many Gabriel-Roiter comeasures. In particular, assume that $\{M_i\}$ is a
family of indecomposable $\MCM$ $R$-modules of bounded $\uh$-length. By Proposition \ref{cor1}, there is a faithful
system of parameters $\xx$ for the family $\{M_i\}$ and so \cite[Corollary 15.11]{lw1} yields that
$\{M_i/{\xx^2M_i}\}$ is a family of indecomposable $R/{\xx^2R}$-modules.
Moreover, according to Lemma \ref{lem9}, the family $\{M_i\}$
is of bounded multiplicity.
Now one may use \cite[Proposition 1.7]{yo}, and conclude that the family
$\{M_i/{\xx^2M_i}\}$ is of bounded length. Consequently,
there are only finitely many Gabriel-Roiter (co)measures for the family $\{M_i/{\xx^2M_i}\}$.\\

{\sc Main Property} (\cite[Main property$^*$]{ri1}). {\em Let $Y_1, \cdots, Y_t, Z$ be indecomposable $\Lambda$-modules of finite
length and assume that there is an epimorphism $g:\oplus_{i=1}^tY_i\lrt Z$.
\begin{itemize}
\item[$(a)$] Then ${\rm min}~{\mu^*}(Y_i)\leq\mu^*(Z)$.
\item[$(b)$] If $\mu^*(Z)={\rm min}~{\mu^*}(Y_i)$, then $g$ splits.

\end{itemize}}
The dual version of the main property (for Gabriel-Roiter measure) has been also appeared in \cite{ri1}.
\end{s}

The next result and the method of its proof, will play an essential role throughout the paper.
\begin{theorem}\label{the2}
Let $\F=\{M_i\}_{i\in I}$ be an infinite set of pairwise non-isomorphic indecomposable $\MCM$ $R$-modules
of bounded $\uh$-length. Then there exists an infinite subset $I'$ of $I$
such that for any $i\in I'$, there is a non-zero $R$-homomorphism
$f_i:M_i\lrt k$ such that for
any $i\neq j\in I'$, any composition map $M_j\lrt M_i\st{f_i}\lrt k$ is zero.
\end{theorem}
\begin{proof}Let us divide the proof into three steps.\\
{\sc step 1}. Since $\F$ is of bounded ${\uh}$-length, by
Proposition \ref{cor1}, there is a faithful system of parameters
$\xx$ for $\F$. In view of Property 2 of \ref{s1},
there are only finitely many Gabriel-Roiter comeasures for {$R/{\xx^2R}$-}modules $M/{\xx^2M}$,
$\mu^*(M/{\xx^2 M})$, where $M\in\F$ {and since $\F$ is infinite,
one may choose an infinite subset $\F^\prime$ of $\mathcal{F}$
consisting of all indecomposable $R$-modules $M_i$
with the same Gabriel-Roiter
comeasure $\mu^*(M_i/{\xx^2M_i})$.}

{\sc step 2}.  Set $I'=\{i\in I\mid M_i\in\F'\}$.
Fix $M_i\in\F^\prime$. We would like to show that any $R$-homomorphism $\oplus_{i\neq j\in I'}M_j^{(\Lambda_j)}\lrt M_i$,
{ where $\Lambda_j$ is a set for each $j$,} is not an epimorphism. Assume on the contrary that there is an epimorphism
$\oplus_{i\neq j\in I'}M_j^{(\Lambda_j)}\lrt M_i$. As $M_i$ is finitely generated, one may find
a finite subset $J$ of $I'$, say $J=\{1, 2, \cdots, s\}$ such that
the $R$-homomorphism $\phi={(\phi_j)_{j=1}^s}:\oplus_{j=1}^sM_j^{n_j}\lrt M_i$ is an epimorphism,
implying that $\bar{\phi}:\oplus_{j=1}^s(M_j/\xx^2M_j)^{n_j}\lrt M_i/\xx^2M_i$
is an epimorphism as well. {In view of \cite[Corollary 15.11]{lw1}, each quotient module is indecomposable. Now}
 since $\mu^*(M_i/\xx^2M_i)=\mu^*(M_j/\xx^2M_j)$
for any $1\leq j\leq s$, by {part $(b)$ of} Main property of \ref{s1}, $\bar{\phi}$ is a split epimorphism
and so {the Krull-Remak-Schmidt theorem} gives rise to the isomorphism
${\bar{\phi}_j}:M_j/\xx^2M_j\lrt M_i/\xx^2M_i$, for some $j\in J$.
{Consequently, by \cite[Lemma 3.3.2]{bh}, $\phi_j:M_j\lrt M_i$ will be an isomorphism.
But this contradicts the hypothesis that modules in $\F$ (and so $\F'$) are non-isomorphic.}

{\sc step 3}. We prove that for any $i\in I'$, there is a non-zero $R$-homomorphism
$f_i:M_i\lrt k$ such that for any $i\neq j\in I'$, any composition map $M_j\lrt M_i\st{f_i}\lrt k$  is zero.
Set $K=\langle\im\phi\rangle{=\Sigma_{\phi}\im\phi}$, where $\phi$ runs over all $R$-homomorphisms
$\phi:\oplus_{i\neq j\in I'}M_j^{(\Lambda_j)}\lrt M_i$.
According to the proof of the previous step, $M_i/K$ is non-zero and so there is a
non-zero homomorphism $g:M_i/K\lrt k$. Therefore, $f_i=g\pi _i:M_i\lrt k$ is
a non-zero $R$-homomorphism, where $\pi_i:M_i\lrt M_i/{K}$ is the natural epimorphism.
Moreover, it is obvious from the construction
of $R$-homomorphisms $f_i^{,}s$ that $M_j\lrt M_i\st{f_i}\lrt k$ is zero, for any $i\neq j\in I'$.
So the proof is completed.
\end{proof}

\section{Balanced big $\CM$ modules with bounded $\CM$-support}
The main theme of this section is to show that every balanced big $\CM$ module
with an $\m$-primary cohomological annihilator and of bounded $\CM$-support, is fully decomposable.

It follows from the definition that balanced big $\CM$ modules
need not be closed under direct summands in general. The next result leads us
to provide a criterion to fix this restriction; see Corollary \ref{cor3}.

\begin{lemma}\label{lem1}Let $M$ be a weak balanced big $\CM$ $R$-module
and $\xx=x_1,\cdots,{x_t}\in\m$ an {$R$-sequence}
 such that $\xx\Ext_R^1(M,-)=0$ over $\Md R$.
If $M/{\xx M}$ is a projective $R/{\xx R}$-module, then $M$ is
projective as an $R$-module. {In particular, if $M/{\xx M}=0$, then $M=0$.}
\end{lemma}
\begin{proof}We prove by induction on {$t$}. Assume that $t=1$.
Since $M/{{x_1} M}$ is a projective $R/{{x_1} R}$-module, we have $\pd_RM/{{x_1} M}\leq 1$
and so $\Ext_R^i(M/{{x_1} M}, -)=0$ over $\Md R$, for any $i\geq 2$.
So by applying the functor $\Hom_{R}(-, N)$, where $N$ is in $\Md R$,
to the short exact sequence of $R$-modules;
$0\lrt M\st{x_1}\lrt M\lrt M/{x_1 M}\lrt 0$, one obtains an exact sequence
$\Ext_R^1(M, N)\st{x_1}\lrt\Ext_R^1(M, N)\lrt 0$.
This means that $\Ext_R^1(M, N)={x_1}\Ext_R^1(M, N)$.
By the hypothesis, the right hand side vanishes,
implying that $\Ext_R^1(M, -)=0$ and then $M$ is projective over $R$.
Now suppose that $t>1$ and the result has been proved for all values smaller than $t$.
Setting $S=R/{\xx_{t-1}R}$, where $\xx_{t-1}=x_1, \cdots,x_{t-1}$,
we have $\pd_SR/{\xx R}\leq 1$, implying that
$\pd_SM/{\xx M}\leq 1$, as well. Considering the following
exact sequence of $S$-modules;
$$0\lrt M/{\xx_{t-1}M}\st{x_t}\lrt M/{\xx_{t-1}M}\lrt M/{\xx M}\lrt 0,$$
one gets an isomorphism $x_t\Ext_S^1(M/{\xx_{t-1}M}, -)=\Ext_S^1(M/{\xx_{t-1}M},-)$ over $\Md S$.
Now one may apply \cite[Lemma 2 (ii) page 140]{ma}, in order to conclude that the isomorphism
$\Ext_S^1(M/{\xx_{t-1}M}, -)\cong\Ext_R^1(M, -)$ holds true over $\Md S$.
On the other hand, by the hypothesis, $x_t\Ext_R^1(M, -)=0$. All of these
facts enable us to deduce that
$\Ext_S^1(M/{\xx_{t-1}M}, -)=0$ over $\Md S$,
meaning that $M/{\xx_{t-1}M}$ is projective over $S$.
Therefore, induction hypothesis would imply that $M$ is indeed projective
over $R$, as desired. {Next assume that $M/{\xx M}=0$. So, $M$ will be a projective $R$-module. Indeed $M$ is a free $R$-module. Now, since $M=\xx M$, one may infer that $M=0$. The proof is completed.} 
\end{proof}

\begin{lemma}\label{lem11}Let $M$ be a weak balanced big $\CM$ $R$-module and
$\xx=x_1, \cdots, x_d$ a weak $M$-regular sequence. If $M=\xx M$,
then for any integer $t>1$, $M=\xx^tM$.
\end{lemma}
\begin{proof}If there is an integer $1\leq i\leq d$ such that $M=x_iM$,
then it is evident that for any integer $t>1$, $M=x_i^tM$ and so, the
equality $M=\xx^tM$ follows. Suppose that, for any $i$, $M\neq x_iM$.
Letting $\xx_{d-1}=x_1,\cdots, x_{d-1}$, the
hypothesis gives rise to the isomorphism; $M/{\xx_{d-1}M}\st{x_d}\lrt M/{\xx_{d-1}M}$.
In particular, the composition map $M/{\xx_{d-1}M}\st{x_d}\lrt M/{\xx_{d-1}M}\st{x_d}\lrt M/{\xx_{d-1}M}$
is again an isomorphism. Indeed, by continuing this for $t$ times, we conclude that
$M/{\xx_{d-1}M}\st{x_d^t}\lrt M/{\xx_{d-1}M}$ is an isomorphism, meaning that
$M=(x_1,\cdots,x_{d-1}, x_d^t)M$. Since
any permutation of $\xx$ is again a weak $M$-regular sequence,
continuing this manner for any $i$, will complete the proof.
\end{proof}

\begin{cor}\label{cor3}
Let $M$ be a balanced big $\CM$ $R$-module with an $\m$-primary cohomological
annihilator. Then any non-zero direct summand of $M$ is balanced big $\CM$.
\end{cor}
\begin{proof}Assume that $M'$ is a non-zero direct summand of $M$. Since $M$ is balanced
big $\CM$, $M'$ is clearly a weak balanced big $\CM$ $R$-module.
Take an arbitrary system of parameters $\xx$ of $R$. If $M'\neq\xx M'$,
then we are done. So assume that $M'=\xx M'$.
As $M$ has an $\m$-primary cohomological annihilator, there is an
integer $t>0$ such that $\xx^t\Ext_R^1(M, -)=0$ and so $\xx^t$ annihilates
the functor $\Ext_R^1(M', -)$, as well. {Since} $M'/{\xx M'}=0$, by Lemma \ref{lem11},
$M'/{\xx^tM'}=0$ and thus Lemma \ref{lem1} ensures that $M'$ is a projective
$R$-module.
Hence the proof is finished.
\end{proof}


\begin{theorem}\label{prop8}
{Let $M$ be a weak balanced big $\CM$ $R$-module with  bounded $\uh$-length on  $\CM$-$\Supp_R(M)$. Let $\xx$ be  a faithful system of parameters  for $\CM$-$\Supp_R(M)$ such that $M/{\xx^2M}$ is non-zero. Then the following hold.
\begin{enumerate}\item There exists an indecomposable $\MCM$ $R$-module $X$ and a non-zero pure monomorphism $\varphi:X\lrt M$. In particular, $\bar{\varphi}:X/{\xx^2X}\lrt M/{\xx^2M}$ is a split monomorphism.
\item If, in addition, $\xx\Ext_R^1(M, N)=0$ for all $\MCM$ modules $N$, then $\varphi$ will be a split monomorphism.
\end{enumerate}}
\end{theorem}
\begin{proof}
(1) Since $M$ is a weak balanced big $\CM$ $R$-module, 
by \cite[Theorem B]{ho} there is a direct system
$\{M_i, \varphi^i_j\}_{i,j\in I}$ of $\MCM$ $R$-modules such that $M=\underrightarrow{\lim}M_i$.
By our assumption,  $\underrightarrow{\lim}M_i/{\xx^2M_i}= M/{\xx^2 M}$ is non-zero.
In what follows, $-\otimes_RR/{\xx^2R}$ for simplicity will be denoted by $\bar{(-)}$.
Thus, we may take an index $j\in I$ and an indecomposable $\MCM$ direct summand
$X_j$ of $M_j$ such that the morphism ${\bar{\varphi}'_j}:\bar{X}_j\lrt \bar{M}$
is non-zero, where ${{\varphi'}_j}={{\varphi}_j}{\mid_{{X_j}}}$ and $\varphi_j:M_j\lrt M$ {is the natural morphism
such that for any $i\leq j$, the equality $\varphi_i=\varphi_j\varphi^i_j$} holds.
Let $k_1 \in I$ be an index with $k_1 >j$, so we have the morphism $\bar{\varphi}_{k_1}^j:\bar{M_j} \lrt\bar{M}_{k_1}$.
In view of the equality $\varphi_j= \varphi_{k_1}\varphi_{k_1}^j$, one may find
an indecomposable $\MCM$ direct summand $X_{k_1}$ of $M_{k_1}$
{such that the composition map
$$\bar{X}_j\st{\bar{\varphi}_{k_1}^j|_{\bar{X}_j}}\lrt\bar{M}_{k_1}\st{{\bar{\pi}}}\lrt\bar{X}_{k_1}
\st{\bar{\varphi}_{k_1}|_{\bar{X}_{k_1}}}\lrt\bar{M}$$
is non-zero, where ${{\pi}} : {M}_{k_1}\lrt {X}_{k_1}$ is the canonical projection.}
We denote the composition map {
 $${X}_j\st{{\varphi}_{k_1}^j|_{{X}_j}}\lrt {M}_{k_1}\st{\pi}\lrt {X}_{k_1}$$ by $\psi^j_{k_1}$.}
Now apply the induction argument to obtain a chain of morphisms of indecomposable finitely generated modules {
$$ \bar{X}_j\st{\bar{\psi}^j_{k_1}}\lrt \bar{X}_{k_1}\st{\bar{\psi}^{k_1}_{k_2}}\lrt \bar{X}_{k_2}\st{\bar{\psi}^{k_2}_{k_3}}\lrt \bar{X}_{k_3}\lrt \cdots,$$}
such that  the compositions  have non-zero images in $\bar{M}$.
This, in particular, means that all $X_i^,$s belong to $\CM$-$\Supp_R(M)$
and so they are of bounded $\uh$-length. Thus by applying Lemmas \ref{lem9} and \ref{lem7}, we infer that there is a non-negative
integer $b$ such that $l_R(\bar{X}_{i})<b$.
Now Harada-Sai Lemma \cite[Lemma 11]{hs},
guarantees the existence of  an index $k_t \in I$ such that for each $k_s>k_t$
the induced morphism {$\bar{\psi}^{k_t}_{k_s}: \bar{X}_{k_t}
\lrt \bar{X}_{k_s}$ needs to be an isomorphism.  So by making use of \cite[Lemma 3.3.2]{bh}, we conclude that ${\psi}^{k_t}_{k_s}: {X}_{k_t}
\lrt {X}_{k_s}$ is an isomorphism, as well. This yields that, for any $k_s> k_t$ the morphism ${\varphi}_{k_s}^{k_t}|_{{X}_{k_t}}: {X}_{k_t} \lrt {M}_{k_s}$ is a split monomorphism.
This, in turn, would imply that ${\varphi}'_{k_t}:{X}_{k_t}\lrt {M}$ is a pure monomorphism.}
As $\bar{X}_{k_t}$ is a finitely generated module over the artinian ring $\bar{R}$, it will be pure injective, inforcing $\bar{\varphi}'_{k_t}$ to be a split monomorphism, giving the desired result.\\ (2) In view of part (1), there is an  indecomposable $\MCM$ module $X$ {and a pure monomorphism $\varphi:X\lrt M$. So $\bar{\varphi}:\bar{X}\lrt \bar{M}$ is a split monomorphism.} Suppose that $g:\bar{M}\lrt \bar{X}$ is an $\bar{R}$-homomorphism with $g{\bar{\varphi}}=id_{\bar{X}}$.
By our assumption, $\xx\Ext_R^1(M, X)=0$ and so a verbatim pursuit of the argument given in the proof of \cite[Proposition 14.9]{lw1}
(see also \cite[Proposition 6.15]{yo}),
yields that there exists an $R$-homomorphism
$h:M\lrt X$ such that $g\otimes R/\xx R=h\otimes R/\xx R$.
Therefore, it is fairly easy to see that $h\varphi\otimes R/\xx R=id_{X/\xx X},$
{and so by \cite[Lemma 3.3.2]{bh}, $h\varphi$ is an isomorphism. Hence $\varphi$ will be a split monomorphism.
The proof now is finished.}
\end{proof}

\begin{remark}\label{cor4}
Let $M$ be as in the above theorem. The proof of Theorem \ref{prop8}, reveals that
for any non-zero element $z\in M/{\xx^2 M}$, there is some indecomposable direct summand $X$ of $M$ such that
$X/{\xx^2 X}$ is a direct summand of $M/{\xx^2 M}$ and $z$ has non-zero component in $X/{\xx^2 X}$,
where $\xx$ is a faithful system of parameters for $X^,$s.
\end{remark}

\begin{s}Let $\Lambda$ be an artin algebra. A result due to Ringel \cite[Theorem 4.2]{ri3}
asserts that an indecomposable $\Lambda$-module $X$ of finite length
with $\mu(X)=\gamma$ is relative $\Sigma$-injective in
$\D(\gamma)$, where $\D(\gamma)$ is the full subcategory
consisting of all $\Lambda$-modules $M$ in which any indecomposable submodule
$M'$ of $M$ of finite length satisfies $\mu(M')\leq\gamma$.
That is to say, any submodule $M'$ of a module $M\in\D(\gamma)$
which is a direct sum of copies of $X$ will be a direct summand of $M$.
By the aid of a counterexample, he realized that the hypothesis $M'$ being
a direct sum of a finite number of non-isomorphic indecomposable modules of finite length
with a fixed Gabriel-Roiter measure $\gamma$, is essential.
Indeed, he showed that there are infinitely many isomorphism classes of
submodules $M_i$ of a module $M\in\D(\gamma)$ such that for any $i$,
$\mu(M_i)=\gamma$, but the embedding $\varphi:\oplus_iM_i\lrt M$ is not split.
The argument given in the proof of the next result reveals that if
we impose the hypothesis that $\mathsf{supp}_{\Lambda}(M)$ is of bounded length, then
$\varphi$ will be split.
\end{s}

\begin{theorem}\label{th7}
{Let $M$ be a weak balanced big $\CM$ $R$-module with  bounded $\uh$-length on  $\CM$-$\Supp_R(M)$. Let $\xx$ be  a faithful system of parameters  for $\CM$-$\Supp_R(M)$. Then the following statements hold.
\begin{enumerate}\item If $M/{\xx^2M}$ is non-zero, then there is a non-zero fully decomposable balanced big $\CM$ module 
$Y$ and a pure monomorphism $\varphi:Y\lrt M$  such that $\bar{\varphi}:Y/{\xx^2Y}\lrt M/{\xx^2M}$ is an isomorphism.
\item If, in addition, $\xx\Ext_R^1(M, -)=0$, then $\varphi$ is an isomorphism.
\end{enumerate}}
\end{theorem}
\begin{proof}(1)
 { By Theorem \ref{prop8}(1), there is a pure monomorphism $i_X:X\lrt M$, where $X$ is an indecomposable $\MCM$ module.} 
 Assume that $\Sigma$ is the set of all pure submodules of $M$ which are direct sums of
indecomposable {$\MCM$ modules}. For any two objects $N, L\in\Sigma$, we write
$N\leq L$ if and only if $N$ is a {pure submodule}  of $L$ and the following diagram
is commutative; {\footnotesize{\[\xymatrix@C-0.5pc@R-.8pc{N\ar[dr]_{i_N} \ar[rr]^{i_{NL}} && L\ar[dl]^{i_L} \\ & M  & }\]}}
where {the inclusion maps are pure monomorphism.}
By Zorn's lemma one may find a pure submodule $Y=\oplus X_i$ of $M$ where each
$X_i$ is an indecomposable $\MCM$ {pure submodule }  of $M$ and $Y$ is maximal with respect
to this property. Take the pure exact sequence of $R$-modules;
$$\eta:0\lrt Y\st{i_Y}\lrt M\st{\upsilon}\lrt K\lrt 0.$$
{So,  the desired result will be achived, if we show that $K/{\xx^2K}=0$.
Assume for the contradiction that this is not the case.}
It is evident that any element of $\CM$-$\Supp_R(K)$ belongs to
$\CM$-$\Supp_R(M)$, and so $\CM$-$\Supp_R(K)$ will be of bounded $\uh$-length
and also $\xx$  is a faithful system of parameters for $\CM$-$\Supp_R(K)$.
Since $\eta$ is a pure exact sequence and the two modules $Y, M$ are balanced big $\CM$, by Lemma \ref{re3}, 
$K$ is a weak balanced big $\CM$ $R$-module. 
 {According to Theorem \ref{prop8}(1), there is a non-zero pure monomorphism $\theta:N\lrt K$, where $N$ is an indecomposable $\MCM$ module.
 Considering the pure exact sequence $\eta$ and the morphism $\theta:N\lrt K$, one may obtain the induced map $\psi:N\lrt M$ such that $\nu\psi=\theta$. As ${\theta}$ is a pure monomorphism, the same will be true for ${\psi}$.}
In particular, we will have
the following commutative diagram;
{\footnotesize{ \[\xymatrix@C-0.8pc@R-.8pc{0\ar[rr]&& Y\ar[rr]^{i_Y}&& M \ar[rr] && K\ar[rr]&& 0, \\ &&&&& N\ar[ur]_{\theta} \ar[ul]^{{\psi}}& & && }\]}} where $\theta$ {and $\psi$ are pure monomorphism. }
 Hence, one may deduce that $Y\oplus N\st{[i_Y~~{\psi}]}\lrt M$ is indeed a
pure monomorphism and $Y\oplus N$ contains $Y$ properly, but this contradicts the maximality of $Y$.
Thus $K/{\xx^2K}=0$ and so  $i_Y\otimes R/{\xx^2R}$ is an isomorphism.
{Now we set $\varphi$ to be $i_Y$, and the desired result is obtained.} 
\\
(2) {By virtue of part (1), the morphism   $\varphi\otimes R/{\xx^2R}:Y/{\xx^2 Y}\lrt M/{\xx^2M}$ is an isomorphism.
 Assume that $\rho:M/{\xx^2M}\lrt Y/{\xx^2Y}$ is the inverse of $\varphi\otimes R/{\xx^2R}$.} 
 {By the hypothesis,}
 $\xx\Ext_R^1(M, Y)=0$, and so we may find a morphism $g:M\lrt Y$
such that $g\otimes_{R} R/\xx R=\rho\otimes_{R} R/\xx R$.
Now we show that $g$ is an isomorphism. Since $\bar{g}:M/{\xx^2M}\lrt Y/{\xx^2Y}$ is an isomorphism,  Lemma \ref{re2} ensures that $g$ is an epimorphism.
 Taking the exact sequence of $R$-modules; $0\lrt L\lrt M\st{g}\lrt Y\lrt 0$,
and using the fact that both modules $M, Y$ have $\m$-primary cohomological annihilators, Remark \ref{re1} yields  that the same is true for the weak balanced big $\CM$ module $L$. {On the other hand, by \cite[Proposition 1.1.5]{bh}, this sequence remains exact after applying the functor $-\otimes_RR/{\xx^2R}$. Consequently,} $L/{\xx^2 L}=0$, {and so,} Lemma \ref{lem1} forces $L$ to be zero, implying that $g$ is an isomorphism. {Hence, we will have the equality $g\varphi\otimes R/{\xx R}=id_Y\otimes R/{\xx R}$, and then, the argument appeared just above, yields that $g\varphi$ is an isomorphism. In particular, $\varphi$ is an isomorphism, becasue $g$ is so.} Thus the proof is completed.
\end{proof}

\begin{s}\label{s3}Anderson and Fuller \cite{af} posed the problem of determining over which
rings does every module has a decomposition $M=\oplus_{i\in I}M_i$ that complements
direct summands in the sense that whenever $K$ is a direct summand
of $M$, $M=K\oplus(\oplus_{j\in J}M_j)$ for some $J\subseteq I$.
This problem has been settled for artin rings of finite
representation type by Tachikawa \cite{t}. The result below indicates that
Tachikawa type theorem satisfies for Cohen-Macaualy rings of finite $\CM$-type.
\end{s}

\begin{theorem}\label{th2} Any balanced
big $\CM$ module with an $\m$-primary cohomological annihilator
and bounded {$\uh$}-length on $\CM$-support, satisfies complements direct summands.
\end{theorem}
\begin{proof} Take a balanced big $\CM$ $R$-module $M$ with an $\m$-primary
cohomological annihilator such that its $\CM$-support is of bounded $\uh$-length.
{By Theorem \ref{th7}, $M$ is fully decomposable and so,
$M=\oplus_{i\in I}M_i$, where each $M_i$ is an indecomposable
finitely generated submodule of $M$. Now assume that $K$ is a direct summand of $M$.}
In view of Corollary \ref{cor3}, any direct summand of $M$ is again balanced big
$\CM$ which has an $\m$-primary cohomological annihilator with bounded
$\uh$-length on $\CM$-support, so another use of Theorem \ref{th7} forces it to
be fully decomposable. {This, in particular, gives rise to another
decomposition of $M$.} Hence the Krull-Schmidt-Azumaya
theorem \cite[pages 331-332]{hs} gives the desired result.
\end{proof}

\begin{cor} Let $R$ be of finite $\CM$-type. Then any balanced big
$\CM$ $R$-module $M$ with an $\m$-primary cohomological annihilator,
is of finite $\CM$-type. In particular,  each balanced
big $\CM$ module with an $\m$-primary cohomological annihilator, satisfies complements direct summands.
\end{cor}

\begin{cor}Let $R$ be an isolated singularity containing its residue field $k$
and let $M$ be a balanced big $\CM$ $R$-module with an $\m$-primary cohomological annihilator such that
$\CM$-$\Supp_R(M)$ is of bounded multiplicity. If $k$ is perfect, then $M$ is fully decomposable.
\end{cor}
\begin{proof}According to \cite[Theorem 14.19]{lw1} (see also \cite[Corolary 2.8]{yo1}),
there is a faithful system of parameters
$\xx$ for the class of all indecomposable $\MCM$ $R$-modules.
In particular, $\xx$ is a faithful system of parameters for $\CM$-$\Supp_{R}(M)$.
So by making use of Lemma \ref{lem88}, $\CM$-$\Supp_R(M)$ is of bounded $\uh$-length.
Now Theorem \ref{th7} gives the desired result.
\end{proof}

\section{$\MCM$ modules of large finite $\mathsf{h}$-length}

This section reveals that any balanced big $\CM$ module with an $\m$-primary
cohomological annihilator is of finite $\CM$-type, whenever $R$ is Gorenstein
or it is an isolated singularity and the class of all indecomposable $\MCM$ $R$-modules
is of bounded $\uh$-length. Our results provide the first Brauer-Thrall
type theorem for rings, concerning the invariant $\uh$-length.

\begin{prop}\label{prop10}Let $\F=\{M_i\}_{i\in I}$ be a set of pairwise non-isomorphic
indecomposable $\MCM$ $R$-modules and let $(f_i:M_i\lrt k)_{i\in I}$ be a family of non-zero
$R$-homomorphisms such that any composition map $M_j\lrt M_i\st{f_i}\lrt k$ with $j\neq i$, is zero.
{Let, for any $i$, $g_i:M_i\lrt G$ be an induced homomorphism by $f_i$ {(i.e. $\alpha g_i=f_i$)}
and set $g=(g_i)_{i\in I}:\oplus_{i\in I}M_i\lrt G$. Assume that $\beta:R^n\lrt G$
is a homomorphism such that $(\oplus_{i\in I}M_i)\oplus R^n\st{[g~~~\beta]}\lrt G$ is an epimorphism.}
If the kernel of {this} epimorphism 
is fully decomposable, then $\F$ is a finite set.
\end{prop}
\begin{proof}Assume on the contrary that $\F$ is an infinite set.
Consider the short exact sequence of $R$-modules,
$0\lrt K\lrt(\oplus_{i\in I}M_i)\oplus R^n\st{[g~~~\beta]}\lrt G\lrt 0$.
By the hypothesis, $K=\oplus_{i\in J}K_i$, where each $K_i$ is an indecomposable
finitely generated module.
Since $K$ is weak balanced big $\CM$, for any $i$, $K_i$ is a $\MCM$ $R$-module.
Therefore, for each $i$, there is an $R$-monomorphism $\epsilon_i:K_i\lrt\omega^{n_i}$
and so we may obtain the following commutative diagram;
\begin{equation}\label{eq1}
{\footnotesize{\begin{CD}
0 @>>> \oplus_{i\in J}K_i @>s>> \ (\oplus_{i\in I}M_i)\oplus R^n @>[g~~~\beta]>> \ G @>>> 0\\
& & @V id VV @V u VV @V \varphi VV & &\\ 0 @>>> \oplus_{i\in J}K_i @>\oplus\epsilon_i >> \oplus_{i\in J}\omega^{n_i} @>\oplus\phi_i>> \oplus_{i\in J}\Omega^{-1}K_i
@>>> 0,\end{CD}}}
\end{equation}
in which the morphism $u$ is induced by the identity map.
As $G$ is finitely generated, the image of $\varphi$ is non-zero only in a
finite number of $\Omega^{-1}K_i^,$s. This, in turn, allows us to decompose
the morphism $\varphi$ into the direct sum of $\varphi':  G\lrt\oplus_{i\in J'}\Omega^{-1}K_i$
and $0\lrt\oplus_{i\in J''}\Omega^{-1}K_i$,
where $J'$ is a finite subset of $J$ and  $J''=J-J'$. {Set, for simplicity, $\epsilon':=\oplus_{i\in J'}\epsilon_i$, $\phi':=\oplus_{i\in J'}\phi_i$, and  we define the morphisms $\epsilon''$ and $\phi'',$ similarly.}
Take the following pull-back diagram;
\begin{equation}
{\footnotesize{\begin{CD}
0 @>>> \oplus_{i\in J'}K_i @>\gamma>> \ M' @>h>> \ G @>>> 0\\
& & @V id VV @V\nu VV @V \varphi' VV & &\\ 0 @>>> \oplus_{i\in J'}K_i @>\epsilon'>>
\oplus_{i\in J'}\omega^{n_i} @>\phi'>> \oplus_{i\in J'}\Omega^{-1}K_i @>>> 0.\end{CD}}}
\end{equation}
{Consider the following commutative diagram; 
{\footnotesize{$$\begin{CD}
0 @>>> \oplus_{i\in J}K_i @>s>> \ (\oplus_{i\in I}M_i)\oplus R^n @>[g~~~\beta]>> \ G @>>> 0\\
& & @V \pi' VV @V u' VV @V \varphi' VV & &\\ 0 @>>> \oplus_{i\in J'}K_i @>\epsilon'>>
\oplus_{i\in J'}\omega^{n_i} @>\phi'>> \oplus_{i\in J'}\Omega^{-1}K_i @>>> 0,\end{CD}$$}}}
{where $\pi'$ is the projection and $u'$ is the composition map
 $(\oplus_{i\in I}M_i)\oplus R^n\st{u}\lrt\oplus_{i\in J}\omega^{n_i}\st{\pi_1}\lrt\oplus_{i\in J'}\omega^{n_i}$, in which $\pi_1$ is the natural projection.}
By using the property of pull-back diagram,
we may find $R$-homomorphisms $\psi:\oplus_{i\in I}M_i\oplus R^n\lrt M'$
{and $t:\oplus_{i\in J}K_i\lrt \oplus_{i\in J'}K_i$} such that
the following diagram; 
\begin{equation}
\label{lk}
{\footnotesize{\begin{CD}
0 @>>> \oplus_{i\in J}K_i @>s>> \oplus_{i\in I}M\oplus R^n @>[g~~\beta]>> \ G @>>> 0\\
& & @V t VV @V\psi VV @V id VV & &\\ 0 @>>> \oplus_{i\in J'}K_i @> \gamma>>
M' @>h>> G @>>> 0,
\end{CD}}}
\end{equation}
{is commutative. Another use of the property of pull-back diagram, gives rise to the equality $u'=\nu\psi$. This, in conjunction with  the commutativity of the left squares in the above two diagrams, leads us to obtain the equality $\epsilon'\pi'=\epsilon' t$. Now $\epsilon'$ being monomorphism, yields that $\pi'=t$. Since $\phi''u=0$, the commutativity of (\ref{eq1}) yields that there exists an $R$-homomorphism
$\theta:(\oplus_{i\in I}M_i)\oplus R^n\lrt\oplus_{i\in J''}K_i$ such that $\epsilon''\theta=u''$, where $u''$ stands for the composition map
 $(\oplus_{i\in I}M_i)\oplus R^n\st{u}\lrt\oplus_{i\in J}\omega^{n_i}\st{\pi''}\lrt\oplus_{i\in J''}\omega^{n_i}$, in which $\pi''$ is the natural projection. Also, another use of the commutativity of (\ref{eq1}) gives rise to the equality $u''s_{\mid\oplus_{i\in J''}K_i}=\epsilon''$, and so $\epsilon''\theta s_{\mid\oplus_{i\in J''}K_i}=\epsilon''$.
As $\epsilon''$ is a monomorphism, $\theta s_{\mid\oplus_{i\in J''}K_i}=id_{\oplus_{i\in J''}K_i}$.}
Therefore, we have  the following commutative diagram;
{\footnotesize{$$\begin{CD}
0 @>>> \oplus_{i\in J}K_i @>s>> \ (\oplus_{i\in I}M_i)\oplus R^n @>[g~~~\beta]>> \ G @>>> 0\\
& & @V {id} VV @V \tiny {\left[\begin{array}{ll} \psi \\
\theta \end{array} \right]} VV @V id VV & &\\ 0
 @>>> (\oplus_{i\in J'}K_i) \oplus (\oplus_{i\in J''}K_i)
@>\gamma\oplus id_{\oplus_{i\in J''} K_i} >> M'\oplus(\oplus_{i\in J''}K_i)  @>[h~~~0]>> G @>>> 0.\end{CD}$$}}
{It should be observed that the commutativity of the right-hand side square
follows from the equality $h\psi=[g~~~\beta]$, however the left-hand side square
is commutative, because of the definition of $\theta$ and the commutativity of the left square in (\ref{lk}), and so, 
 $\tiny {\left[\begin{array}{ll} \psi \\ \theta \end{array} \right]}$
will be an isomorphism with inverse $\eta$. }In particular, one obtains the next commutative square;
{\footnotesize{$$\begin{CD}
\ (\oplus_{i\in I}M_i)\oplus R^n @>[f~~~\alpha\beta]>> \ k \\
@V \tiny {\left[\begin{array}{ll} \psi \\ \theta \end{array} \right]} VV @V id VV & &\\
 M'\oplus(\oplus_{i\in J''}K_i)  @>[\alpha h~~~0]>> k,\end{CD}$$}}
{where $f=(f_i):\oplus_{i\in I}M_i\lrt k$} {and $\alpha:G\lrt k$ is the right minimal $\MCM$-approximation of $k$.}
As $M'$ is finitely generated, 
{there are only} finitely many $M_i^,$s; say
$\{M_{i_1},\cdots, M_{i_t}\},$ {such that under $\eta$, $M'$ may have non-zero image in $\{M_{i_1}, \cdots, M_{i_t}, R^n\}$. Since $\F$ is assumed to be infinite, one may}
take a non-projective indecomposable module $M_s$ in $\F$ such that $s\notin\{i_1, \cdots, i_t\}$.
As $f_s:M_s\lrt k$ is non-zero, there is an element $x\in M_{s}$ such that $f_s(x)\neq 0$
and so the image of $x$ under the composition map
$M_{s}\st{i}\lrt\oplus_{i\in I}M_i\oplus R^n\st{[f~~~\alpha\beta]}\lrt k$
will be non-zero, {where $i$ is the injection map}. Thus the commutativity of the above square
enables us to conclude that the image of $x$, say $x'$, under the morphism
$$M_{s}\st{i}\lrt(\oplus_{i\in I}M_i)\oplus R^n\st{\tiny {\left[\begin{array}{ll} \psi \\ \theta \end{array} \right]}}\lrt
M'\oplus(\oplus_{i\in J''}K_i)\st{\pi}\lrt M'$$
is non-zero. Consequently, $(x', 0)$ is a non-zero element of
$M'\oplus(\oplus_{i\in J''}K_i)$ and in particular, $\alpha h(x')$ is non-zero in $k$, as well. {Therefore,} the composition map;{
$$M_{s}\st{i}\lrt\oplus_{i\in I}M_i\oplus R^n\st{\psi}\lrt M'\st{\eta_{\mid M'}}\lrt \oplus_{j=1}^{t}M_{i_j}\oplus R^n\st{f'}\lrt k,$$ is non-zero, where $f'=f_{\mid {\oplus_{j=1}^{t}M_{i_j}\oplus R^n}}$}. {On the other hand},
the construction of morphisms $f_i^,$s, indicates that the composition map  $$M_{s}\st{\psi i}\lrt M'\st{\eta'}\lrt \oplus_{j=1}^{t}M_{i_j}\lrt\oplus_{i\in I}M_i\st{f}\lrt k,$$ is zero. {Here $\eta'$ stands for the composition map $M'\st{\eta_{\mid M'}}\lrt\oplus_{j=1}^{t}M_{i_j}\oplus R^n\lrt\oplus_{j=1}^{t}M_{i_j}$.}
Consequently, {the composition map $M_s\st{\psi i}\lrt M'\st{\eta''}\lrt R^n\st{\alpha\beta}\lrt k$ will be non-zero, where $\eta''$ is the composition map $M'\st{\eta_{\mid M'}}\lrt\oplus_{j=1}^{t}M_{i_j}\oplus R^n\lrt R^n$. This implies that $M_s$ is isomorphic to $R$, because $M_s$ is indecomposable, which contradicts the choice of $M_s$. Hence $\F$ will be a finite set. The proof then is completed.}
\end{proof}

Now, we are in a position to state the main theorem of this
section, which provides the local version of the first Brauer-Thrall conjecture, i.e.
for modules instead of the base ring. 

\begin{theorem}\label{th1}Let $R$ be an isolated singularity
and let $M$ be a balanced big $\CM$ $R$-module having an $\m$-primary cohomological
 annihilator. If $M$ is not of finite $\CM$-type, then there are indecomposable
$\MCM$ $R$-modules of arbitrary large (finite) $\uh$-length.
\end{theorem}
\begin{proof}
Suppose on the contrary that the class of all indecomposable
$\MCM$ $R$-modules, is of bounded $\uh$-length.
Since the same will be true for $\CM$-$\Supp_R(M)$,
in view of Theorem \ref{th7}, $M$ is fully decomposable. So, we may write
$M=\oplus_{i\in I}M_i^{(t_i)}$, where each $M_i$ is an indecomposable $\MCM$ $R$-module.
As $M$ is not of finite $\CM$-type, $\F=\{M_i\}_{i\in I}$ is an infinite set
of pairwise non-isomorphic indecomposable $\MCM$ $R$-modules.
In addition, $\F$ is of bounded $\uh$-length, because $\CM$-$\Supp_R(M)$ is so.
According to Theorem \ref{the2}, there is an infinite
subset $I'$ of $I$ in which for any $i\in I'$, there exists
a non-zero $R$-homomorphism $f_i:M_i\lrt k$ such that for each
$j\neq i\in I'$, any composition map $M_j\lrt M_i\st{f_i}\lrt k$ is zero.
As $\alpha:G\lrt k$ is a right minimal $\MCM$-approximation, for any $i\in I'$,
one may find an $R$-homomorphism $g_i:M_i\lrt G$ such that $\alpha g_i=f_i$.
Set $g=(g_i)_{i\in I'}:\oplus_{i\in I'}M_i\lrt G$.
Consider the exact sequence of $R$-modules;
$0\lrt K\st{\theta}\lrt (\oplus_{i\in I'}M_i)\oplus R^n\st{[g~~~\beta]}\lrt G\lrt 0$.
We claim that $K$ is a balanced big $\CM$ module. To this end,
suppose that $\yy$ is an arbitrary system of parameters of $R$.
Evidently, $\yy$ is a weak $K$-regular sequence, because it is
regular sequence for both modules $(\oplus_{i\in I'}M_i)\oplus R^n$ and $G$.
In addition, $K\neq{\yy K}$. Indeed, if this is not the case,
we will obtain an isomorphism of $R$-modules;
$(\oplus_{i\in I'}M_i)\oplus R^n/{\yy((\oplus_{i\in I'}M_i)\oplus R^n)}\cong
G/{\yy G}$, and this would be contradiction, because $G/{\yy G}$ is finitely generated
whereas $(\oplus_{i\in I'}M_i)\oplus R^n/{\yy((\oplus_{i\in I'}M_i)\oplus R^n)}$ is not so, and thus
the claim follows. Next $M$ having an $\m$-primary cohomological annihilator, yields
that $\m^t\Ext_R^1((\oplus_{i\in I'}M_i)\oplus R^n,-)=0$ for some integer $t>0$.
On the other hand, as $G$ is locally free on the punctured spectrum,
there is an integer $t'>0$ such that $\m^{t'}\Ext_R^2(G, -)=0$.
Therefore $\m^{t+t'}\Ext_R^1(K, -)=0$, meaning the balanced big $\CM$ $R$-module
$K$ has an $\m$-primary cohomological annihilator. Hence, another use of Theorem \ref{th7}
yields that $K$ is fully decomposable. Namely, $K=\oplus_{i\in J}K_i$, where each $K_i$
is an indecomposable finitely generated $R$-module.
Now, Proposition \ref{prop10} forces $I'$ to be finite, which contradicts
with the fact that $I'$ is infinite. The proof then is completed.
\end{proof}

Here we include several corollaries of Theorem \ref{th1}.
First of all, this theorem enables us to prove the first Brauer-Thrall type theorem
for complete Cohen-Macaulay local rings with considering the invariant $\uh$-length.
It is worth noting that, this conjecture fails in general by the aid of
counterexamples of Dieterich \cite{di} and Leuschke and
Wiegand \cite{lw2}, dealing with multiplicity.
\begin{cor}\label{cor2}Let the category of all indecomposable $\MCM$ $R$-modules
be of bounded $\uh$-length. Then $R$ is of finite $\CM$-type.
\end{cor}
\begin{proof}
{By the hypothesis, any indecomposable $\MCM$ $R$-module $X$ has finite
$\uh$-length and so $\uh(X)$ is an artinian $R$-module. Consequently,
for any $\MCM$ $R$-module $M$, $\uh(M)$ will be also an artinian module, {implying that $\uh(M)_{\p}=0$
for all nonmaximal prime ideals $\p$ of $R$. In particular, we have $(\underline{\Hom}_R(M, M))_{\p}\cong\underline{\Hom}_{R_{\p}}(M_{\p}, M_{\p})=0,$ and so $M_{\p}$ is a free}
$R_{\p}$-module, for all nonmaximal prime ideals $\p$ of $R$ and so
by \cite[Lemma 3.3]{yo}, $R$ is an isolated singularity.}
Now, Suppose for the contradiction
that there is an infinite set $\{M_i\}_{i\in I}$ of pairwise non-isomorphic
indecomposable $\MCM$ $R$-modules. So $M=\oplus_{i\in I}M_i$
is not of finite $\CM$-type. On the other hand, in view of Proposition
\ref{cor1}, there is a system of parameters $\xx$ of $R$ such that $\xx\Ext_R^1(M_i,-)=0$
for any $i\in I$. Consequently, $\xx\Ext_R^1(\oplus_{i\in I}M_i, -)\cong\prod_{i\in I}\xx\Ext_R^1(M_i, -)=0$,
meaning that the balanced big $\CM$ module $M$ has an $\m$-primary cohomological
annihilator. Therefore, by virtue of Theorem \ref{th1},
there exist indecomposable $\MCM$ $R$-modules of arbitrary large
(finite) $\uh$-length, which is a contradiction.
The proof hence is completed.
\end{proof}

The above corollary leads us to deduce a result of
Dieterich \cite{di}, Leuschke and Wiegand \cite{lw2} and Yoshino \cite{yo1}.
Indeed we have the next result.

\begin{cor}Let $(R, \m, k)$ be a complete equicharacteristic Cohen-Macaulay local ring
with algebraically closed residue field $k$. Then $R$ is of
finite $\CM$-type if and only if $R$ is an isolated singularity
and there is a bound on the multiplicities of the indecomposable $\MCM$ $R$-modules.
\end{cor}
\begin{proof}{Since the `only if' part is evident, we prove only the `if' part.}
To do this, according to Corollary \ref{cor2},
it suffices to show that the category of all indecomposable $\MCM$ $R$-modules is of bounded
$\uh$-length. By \cite[Theorem 14.19]{lw1}, $R$ admits a faithful system of parameters $\xx$.
Moreover, by the hypothesis, there is an integer $b>0$ such that $e(M)<b$ for any indecomposable
$\MCM$ $R$-module $M$. Now Lemma \ref{lem88} finishes the proof.
\end{proof}

The result below, can be proved similarly to the above corollary.
\begin{cor}Let $R$ be a $d$-dimensional complete Cohen-Macaulay
local ring containing the residue field that is perfect.
Let $M$ be a balanced big $\CM$ $R$-module having an $\m$-primary
cohomological annihilator which is not of finite $\CM$-type. Then there are indecomposable
$\MCM$ $R$-modules of arbitrarily large multiplicity.
\end{cor}

In the remainder of this section, we want to show that over Gorenstein local rings
of finite $\CM$-type, in Theorems \ref{th7} and \ref{th2}, the hypothesis $M$ having an $\m$-primary cohomological
annihilator, is redundant; see Theorem \ref{th8} and Corollary \ref{cor5}.

\begin{lemma}\label{lem12}Let $N$ be a $\MCM$ $R$-module and
$\xx$ a system of parameters of $R$ such that $\xx\Ext_{R}^1(N^*, -)=0$. Then $\xx\Ext_R^1(M, N)=0$
for any module $M\in\X_{\omega}$, where $N^*=\Hom_R(N, \omega)$.
\end{lemma}
\begin{proof}Since $\Ext_R^i(N^*, \omega)=0$, applying the functor $\Hom_R(-, \omega)$
to a free resolution ${\mathbf P_{\bullet}}:\cdots\lrt P_1\lrt P_0\lrt N^*\lrt 0$ of $N^*$
gives rise to the exact sequence of $R$-modules;
$0\lrt N\lrt\Hom_R(P_0, \omega)\lrt\Hom_R(P_1, \omega)\lrt\cdots.$ Thus, for a given object
$M\in\X_{\omega}$, we have the following isomorphisms;
\[\begin{array}{lllll}
\Ext_R^1(M, N) & \cong \HT^1(\Hom_R(M, \Hom_R({\mathbf{P}_{\bullet}, \omega})))\\
& \cong \HT^1(\Hom_R(\mathbf{P}_{\bullet}, M^*))\\
& \cong \Ext_R^1(N^*, M^*),
\end{array}\]
giving the desired result.
\end{proof}

\begin{lemma}\label{lem10}Let $R$ be of finite $\CM$-type
and $\{N_i\}_{i\in I}$ a family of $\MCM$ $R$-modules. Then
there is an integer $t>0$ such that $\m^t\Ext_R^1(M,\oplus_{i\in I}N_i)=0$,
for any module $M\in\X_{\omega}$.
\end{lemma}
\begin{proof}Assume that $\{ X_1, X_2, \cdots, X_t\}$ is the set of all pairwise non-isomorphic
indecomposable $\MCM$ $R$-modules. Take cardinal numbers $s_1, s_2, \cdots, s_t$
such that $\oplus_{i\in I}N_i=\oplus_{i=1}^tX_i^{(s_i)}$ and  assume that $s$ is a non-negative integer
with  $\m^s\uh(X_i^*)=0$, for any $1\leq i\leq t$.
Suppose that for each $i$, $\mathbf{P}_{X_i^*}$ is a projective resolution of $X_i^*$.
So considering an arbitrary $R$-module $M$ in $\X_{\omega}$, analogues to the proof of Lemma \ref{lem12},
we have the following isomorphisms;
\[\begin{array}{lllll}
\Ext_R^1(M, \oplus_{i\in I}N_i) & \cong \HT^1(\Hom_R(M, \oplus_{i=1}^t\Hom_R(\mathbf{P}_{X_i^*},\omega)^{(s_i)}))
\\ & \cong \oplus_{i=1}^t\HT^1(\Hom_R(M, \Hom_R(\mathbf{P}_{X_i^*}, \omega^{(s_i)})))\\
& \cong \oplus_{i=1}^t\HT^1(\Hom_R(\mathbf{P}_{X_i^*}, \Hom_R(M, \omega^{(s_i)})))\\
& \cong\oplus_{i=1}^t \Ext_R^1(X_i^*, \Hom_R(M, \omega^{(s_i)})),
\end{array}\]
giving the desired result.
\end{proof}

\begin{lemma}\label{lem14}Let $0\lrt M\lrt F\lrt L\lrt 0$
be an exact sequence of $R$-modules such that $F$ is free.
 If $L$ admits a non-projective indecomposable $\MCM$
direct summand, then $M$ has an indecomposable $\MCM$ direct summand.
\end{lemma}
\begin{proof}Assume that $L'$ is a non-projective indecomposable $\MCM$ direct summand of $L$.
Consider the following commutative diagram with exact rows;
\begin{equation}\label{eq11}
{\footnotesize{\begin{CD}
0 @>>> K @>>> \ R^n @>>> \ L' @>>> 0\\
& & @V VV @VVV @V i VV & &\\ 0 @>>> M @> >>F @> >> \ L
@>>> 0\\ & & @V VV @V VV @V \pi VV & & \\ 0 @>>> K @> >>R^n @> f>> \ L'
@>>> 0,\end{CD}}}
\end{equation} where $R^n\lrt L'$ is a projective cover.
Since $L'$ is non-projective, the finitely generated $R$-module $K$ is non-zero.
Now using the fact that the right column is identity
and the middle one is an isomorphism, we infer that the left column will be an isomorphism.
This means that the $\MCM$ $R$-module $K$ is a direct summand of $M$, as required.
\end{proof}

{Recall that a commutative noetherian local ring $R$ is said to be {\em Gorenstein}, if it has finite self-injective dimension.  A (not necessarily finitely generated) module $M$ over a Gorenstein ring $R$ is called {\it Gorenstein projective}, whenever it admits a right resolution of projective modules, i.e. $M\in\X_R$. For the basic properties of these modules, we refer the reader to \cite{ej}.}
In the setting of artinian rings, the result below is \cite[Corollary 6]{mo}.
\begin{prop}\label{lem13}Let $R$ be a Gorenstein ring of finite $\CM$-type.
Then any non-zero Gorenstein projective $R$-module has an indecomposable
$\MCM$ direct summand.
\end{prop}
\begin{proof}
Assume that $M$ is a non-zero Gorenstein projective $R$-module and
$\xx$ is a faithful system of parameters
for the class of $\MCM$ $R$-modules, which exists by Proposition \ref{cor1}.
{Moreover, Lemma \ref{lem10} allows us to further assume that $\xx\Ext_R^1(N, \oplus_{i\in I}Y_i)=0$, for any  Gorenstein projective
module $N$ and any family of $\MCM$ modules $\{Y_i\}_{i\in I}$.}
We prove the result in two steps.\\
Step 1: We show that if $M/{\xx^2M}\neq 0$, then $M$ admits an
indecomposable $\MCM$ direct summand.
{Since $M$ is a Gorenstein projective $R$-module and so it belongs to $\X_R$, one may infer that $M$ is a weak balanced big $\CM$ module. So
in view of Theorem \ref{prop8}(1), there is a non-zero pure monomorphism $\varphi:X\lrt M$, where $X$ is an indecomposable $\MCM$ module.
As $X^*=\Hom_R(X, R)$ is a $\MCM$  $R$-module, $\xx\Ext_R^1(X^*, -)=0$. Now, by applying Lemma \ref{lem12}, we get that $\xx\Ext_R^1(M, X)=0$. Hence the argument given in the proof of Theorem \ref{prop8}(2), reveals that
$\varphi$ is a split monomorphism.}\\
Step 2: We prove that $M/{\xx^2M}\neq 0$.
Suppose for the contradiction that $M/{\xx^2M}=0$.
Take a short exact sequence of $R$-modules;  $0\lrt M\lrt F\lrt L\lrt 0$,
in which $F$ is free and $L$ is Gorenstein projective.
As $F/{\xx^2F}\neq 0$, we conclude that the same will be true for $L/{\xx^2L}$.
{So, by Theorem \ref{th7}(1), there is a pure monomorphism $\varphi:Y=\oplus X_i\lrt L$,
where each $X_i$ is an indecomposable $\MCM$ module, such that $\bar{\varphi}:Y/{\xx^2Y}\lrt L/{\xx^2L}$
is an isomorphism. It is evident that, $X_i\lrt Y\lrt L$, for any $i$, is also pure monomorphism.
So, the argument given in the proof of step (1), indicates that these pure monomorphisms  are split.}
By virtue of Lemma \ref{lem14}, we can assume that any indecomposable $\MCM$ direct summand of $L$
is projective. {Consequently, each $X_i$, and then $Y,$ will be projective $R$-modules.}
 {Assume that $\rho:L/{\xx^2L}\lrt Y/{\xx^2Y}$ is the inverse of 
 {$\bar{\varphi}$} .}
Since $\xx\Ext_R^1(L, Y)=0$,  one may find a morphism $g:L\lrt Y$
such that $g\otimes_{R} R/\xx R=\rho\otimes_{R} R/\xx R$.
By Lemma \ref{re2},
$g$ is an epimorphism. Take a short exact
sequence of $R$-modules, $0\lrt T\lrt L\st{g}\lrt Y\lrt 0$. 
As $g\otimes R/{\xx R}$ is an isomorphism {and by \cite[Proposition 1.1.5]{bh},  this sequence remains exact after applying the functor $-\otimes_R R/{\xx R}$,}
we deduce that $T/{\xx T}=0$ and so, the same is true for $T/{\xx^2T}$,
 thanks to Lemma \ref{lem11}. Consider the following pull-back diagram;
{\footnotesize{\[\xymatrix@C-0.5pc@R-.8pc{&&0\ar[d]&0 \ar[d] & &&\\0\ar[r]& M \ar[r] \ar@{=}[d] & F' \ar[r] \ar[d] & T\ar[r]
  \ar[r] \ar[d] & 0& \\ 0 \ar[r]  & M
\ar[r] & F \ar[r] \ar[d] & L \ar[d] \ar[r]  & 0&\\
  && Y \ar@{=}[r] \ar[d] & Y
 \ar[d] &&&\\ &&0& 0& &&\\
}\]}}
Since {$Y$ is projective, the short exact sequence $0\lrt F'\lrt F\lrt Y\lrt 0$} is split, {and then,} $F'$ will be a projective $R$-module.
On the other hand, as $M/{\xx^2M}=0=T/{\xx^2T}$, one may infer that $F'/{\xx^2F'}=0$,
implying that $F'=0$ and so the same is true for {$M$, which is a contradiction.The proof now is finished.}
\end{proof}

As a direct consequence of Proposition \ref{lem13}, we include the following result.
\begin{cor}\label{cor7}Let $R$ be a Gorenstein ring of finite $\CM$-type.
Then any Gorenstein projective $R$-module is balanced big $\CM$.
\end{cor}

\begin{cor}\label{cor8}Let $R$ be a Gorenstein ring and $M$ a Gorenstein projective $R$-module
such that $\CM$-$\Supp_R(M)$ is of bounded $\uh$-length. If $M$ is balanced big $\CM$,
then $M$ has an indecomposable $\MCM$ direct summand.
\end{cor}
\begin{proof}{According to Proposition \ref{cor1}, we may find a faithful system of parameters $\xx$ for $\CM$-$\Supp_R(M)$ such that $\xx\underline{\Hom}_R(X, X)=0$ for any object $X$ in $\CM$-$\Supp_R(M)$. Set $\X=\{X^*\mid X\in \CM$-$\Supp_R(M) \}$, where $X^*=\Hom_{R}(X, R)$.Take an arbitrary object $X$ of $\CM$-$\Supp_R(M)$. Since $(-)^*$ is a duality on the category of finitely generated Gorenstein projective modules, and in particular, on its stable category modulo projectives, we get the isomorphism $\underline{\Hom}_R(X, X)\cong\underline{\Hom}_R(X^*, X^*)$, and then, $\xx\underline{\Hom}_R(X^*, X^*)=0$. Now one may apply a theorem of Hilton-Rees \cite{hr} and deduce that $\xx\Ext_R^1(X^*, -)=0$. This, indeed, means that $\xx$ is a faithful system of parameters for $\CM$-$\Supp_R(M)\cup\X$.}
Now, repeating the proof of step 1 of Proposition \ref{lem13}, will give the desired result.
\end{proof}

\begin{theorem}\label{th8}
Let $R$ be a complete Gorenstein local ring. If $R$ is of finite $\CM$-type,
then  every Gorenstein projective $R$-module is fully decomposable. In particular,
any Gorenstein projective $R$-module satisfies complements direct summands.
\end{theorem}
\begin{proof}Take an arbitrary Gorenstein projective $R$-module $M$.
By Proposition \ref{cor1}  
there is a faithful system of parameters $\xx$ for $\MCM$ $R$-modules.
{According to Corollary \ref{cor7}, $M$ is a balanced big $\CM$ module, and inparticular,
$M/{\xx^2M}$ is non-zero. Now Theorem \ref{th7}(1) guarantees the existence of 
a pure monomorphism $\varphi:Y=\oplus_{i\in I}X_i\lrt M$, where each $X_i$
is indecomposable $\MCM$,  such that $\bar{\varphi}:Y/{\xx^2Y}\lrt M/{\xx^2M}$ is an isomorphism.}
 { Assume that $\rho:M/{\xx^2 M}\lrt Y/{\xx^2 Y}$
is the inverse of $\bar{\varphi}.$
In view of Lemma \ref{lem10}, $\xx\Ext_R^1(M, \oplus_{i\in I} X_{i})=0$, and so,  one may find an
$R$-homomorphism $g:M\lrt Y$ such that $\rho\otimes_RR/{\xx R}=g\otimes_RR/{\xx R}$.
By making use of Lemma \ref{re2},  we infer that $g$ is an epimorphism. Take a short exact sequence
$0\lrt L\lrt M\st{g}\lrt Y\lrt 0$. As $M, Y$ are Gorenstein projective, so is $L$.
In particular, $L$ is balanced big $\CM$. {Since by \cite[Proposition 1.1.5]{bh}, the functor $-\otimes_R R/{\xx^2 R}$ leaves this sequence exact and $g\otimes R/{\xx^2R}$ is an isomorphism, we infer that} 
 $L/{\xx^2L}=0$. {This, in turn, }
would imply that $L=0$, meaning that $g$ is an isomorphism.}
So we are done.
\end{proof}

\begin{cor}\label{cor5}Let $R$ be of finite $\CM$-type. Then a given module $M$ is
a direct sum of $\MCM$ modules if and only if $M\in\X_{\omega}$ and any non-zero direct
summand of $M$ is balanced big $\CM$.
\end{cor}
\begin{proof}By making  use of Lemma \ref{lem10} in  the proof of Theorem \ref{th7},
one can deduce the `if' part of the result. For the `only if' part,
assume that $M'$ is an arbitrary
non-zero direct summand of $M$. We would have nothing to prove, whenever $M'$ is projective.
So assume that $M'$ is not projective. As {$M$ is fully decomposable, it belongs to $\X_{\omega}$, and then} one may easily infer that
$M'$ is weak balanced big $\CM$. Take an arbitrary system of parameters $\xx$ of $R$.
Since $R$ is of finite $\CM$-type, it will admit a faithful system of parameters,
and so we may assume further that $\xx$ is also a faithful system of parameters for
$\MCM$ modules.
On the other hand, the assumption $M$ being fully decomposable leads us to
deduce that $\xx\Ext_R^1(M', -)=0$, because the same is true for $M$.
Now by making use of Lemma \ref{lem1}, we obtain that $M'$ is a balanced
big $\CM$ module, as needed.
\end{proof}

\section{Representation properties of balanced big $\CM$ modules}

The aim of this section is to show that the representation-theoretic properties of balanced big $\CM$ modules
have important consequences for the structural shape of the ring.
It will turn out that any balanced big $\CM$ $R$-module which belongs to $\X_{\omega}$,
being fully decomposable forces $R$ to be an isolated singularity.
Moreover, it is proved that $R$ is of finite $\CM$-type
if and only if it is an isolated singularity and the category of all fully decomposable
modules is closed under kernels of epimorphisms. First we state a lemma.

\begin{lemma}\label{lem8}Let $A$ be a noetherian ring
and let $\{M_i, \varphi^i_j\}_{i,j\in I}$ be a direct
system of indecomposable finitely generated
$A$-modules, over a totally ordered set $I$. If $\underrightarrow{\lim}M_i$ is a direct sum of finitely
generated modules, then $\underrightarrow{\lim}M_i=0$ or there exists an index $t\in I$
such that for any $i\geq t$,  $\varphi^t_i$ is an isomorphism.
\end{lemma}
\begin{proof}If $\underrightarrow{\lim}M_i=0$, then  there is nothing to prove.
So assume that $\underrightarrow{\lim}M_i$ is non-zero.
{By our asssumption, $\underrightarrow{\lim}M_i=\oplus_{j\in J}C_j$,
where each $C_j$ is a finitely generated $A$-module. Since }
by \cite[Corollary 1.2.7]{gt}, $\eta:0\lrt L\lrt\oplus_{i\in I}M_i\st{\varphi=(\varphi_i)_{i\in I}}\lrt \underrightarrow{\lim}M_i\lrt 0$
is a pure exact sequence, the functor $\Hom_A(C_j, -)$ leaves the previous sequence exact, for any $j$,
implying that $\Hom_A(\oplus_{j\in J}C_j, -)$ leaves this sequence exact as well.
This indeed means that $\eta$ is split. Take an $A$-homomorphism
$\psi=(\psi_i)_{i\in I}:\underrightarrow{\lim}M_i\lrt\oplus_{i\in I}M_i$
with $\varphi\psi=\Sigma_{i\in I}\varphi_i\psi_i=id_{\underrightarrow{\lim}M_i}$.
Take a non-zero finitely generated module $C_j$
{such that, under $\psi$, it has non-zero image in only finitely many of $M_{i}^,$s, say $M_1, \cdots, M_t$}.
So, we may define an $A$-homomorphism $\psi'_t:C_j\st{i}\lrt\underrightarrow{\lim}M_i\st{\psi}\lrt\oplus_{i\in I}M_i\st{\beta_t}\lrt M_t,$
where $i$ is injection and {$\beta_t((y_i)_{i\in I})=\Sigma_{i=1}^{t}\varphi_t^i(y_i)$,
for any $(y_i)_{i\in I}\in\oplus_{i\in I}M_i$}.
Considering an $A$-homomorphism $\varphi'_t:M_t\st{\varphi_t}\lrt
\underrightarrow{\lim}M_i\st{\rho}\lrt C_j$,
where $\rho$ is the projection map, we have
$\varphi'_t\psi'_t=\varphi'(\beta_t\psi i)=\rho(\varphi_t\beta_t\psi)i$.
Now, by making use of the following commutative diagram;
{\footnotesize{ \[\xymatrix{\oplus_{i=1}^tM_i\ar[dr]_{\beta_t} \ar[rr]^{\varphi_{| \oplus_{i=1}^t M{_i}}} && \underrightarrow{\lim}M_i \\ & M_t \ar[ur]_{\varphi_t} & }\]}} and the fact that $\varphi\psi=id_{\underrightarrow{\lim}M_i}$, one may infer that $\varphi'_t\psi'_t=id_{C_j}$.
Therefore, $\psi'_t$ is a split monomorphism and so
it will be an isomorphism, because $M_t$ is indecomposable.
Since for any $s\geq t$ we have the following commutative diagram
of $A$-modules;
{\footnotesize{\[\xymatrix{M_t\ar[dr]_{\varphi^t_s} \ar[rr]^{\varphi_t} && \underrightarrow{\lim}M_i \\ & M_s \ar[ur]_{\varphi_s} & }\]}}
One may have the equalities;
$id_{C_j}=\varphi'_t\psi'_t=\rho\varphi_t\psi'_t=(\rho\varphi_s)(\varphi^t_s\psi'_t).$
Thus we have that $\varphi^t_s\psi'_t:C_j\lrt M_t\lrt M_s$
is a split monomorphism, for any $s\geq t$. Now $C_j$ and $M_s$ being indecomposable, forces
$\varphi^t_s\psi'_t$ to be an isomorphism. This, in turn,
implies that $\varphi^t_s$ is indeed an isomorphism. The proof then is completed.
\end{proof}
\begin{lemma}\label{prop6}Let $\{M_i, \varphi^i_j\}_{i,j\in\N_0}$ be a direct system of
$\MCM$ $R$-modules such that for any $i\leq j$, $\varphi^i_j:M_i\lrt M_j$ is a monomorphism
with $\cok\varphi^i_j$ is $\MCM$. Then $\underrightarrow{\lim}M_i$ belongs to $\X_{\omega}$.
\end{lemma}
\begin{proof}First one should note that, { as we have mentioned in \ref{s4}(1), each $\MCM$ $R$-module belongs to $\X_{\omega}$.}
{For any $R$-module $M_i$, we inductively construct a right resolution by modules in $\add\omega$
forming a direct system. {Set $M_{-1}=0$}.  For $i=0$, take an exact sequence of $R$-modules;
$0\lrt M_0\lrt\omega_{0_0}\lrt\omega_{0_1}\lrt\cdots$, where for any $j\geq 0$,
$\omega_{0_j}\in\add\omega$.
Now assume that $i\geq 0$ and we have constructed such a resolution for $M_i$ with the following commutative diagram of $R$-modules;
{\footnotesize {{\[\xymatrix{0 \ar[r] & M_{i-1} \ar[r] \ar[d]^{\varphi^{i-1}_i} &
 \omega_{(i-1)_0} \ar[r]\ar[d] & \omega_{(i-1)_1} \ar[r] \ar[d] & \omega_{(i-1)_2}
\ar[r] \ar[d] & \cdots
\\ 0 \ar[r] & M_i \ar[r] & \omega_{i_0} \ar[r] & \omega_{i_1} \ar[r] & \omega_{i_2} \ar[r] & \cdots, }\]}}}
in which, all columns, except the left-hand side, are split monomorphism.
Now we construct the diagram for the case $i+1$.
By our assumption, $\cok\varphi^i_{i+1}$ is a $\MCM$ $R$-module and so,
there is an exact sequence of $R$-modules;
$0\lrt\cok\varphi^i_{i+1}\lrt\omega'_{i_0}\lrt\omega'_{i_1}\lrt\cdots,$ such that
each $\omega'_{i_j}$ lies in $\add\omega$.
Since the functor $\Hom_R(-, \omega)$
leaves any short exact sequence in $\MCM$ modules, exact,
one may construct the following commutative
diagram of $R$-modules with exact rows and columns;
{\footnotesize{\[\xymatrix@C-0.5pc@R-.8pc{&0\ar[d]&0 \ar[d] &0 \ar[d] &&\\0\ar[r]& M_i
 \ar[r] \ar[d]^{\varphi^i_{i+1}} & \omega_{i_0}  \ar[d]\ar[r] & \omega_{i_1}  \ar[d]\ar[r]&
\cdots& \\ 0 \ar[r]  & M_{i+1} \ar[r] \ar[d] & \omega_{{(i+1)}_{0}} \ar[d] \ar[r] & \omega_{{(i+1)}_{1}} \ar[d] \ar[r]&\cdots&\\
0\ar[r] & \cok\varphi^i_{i+1}\ar[r] \ar[d] & \omega'_{i_0}\ar[r] \ar[d]
 &\omega'_{i_1} \ar[r] \ar[d]&\cdots,&\\ &0& 0& 0& &\\}\]}}
such that any column, except the left-hand side, is split.}
Hence applying the direct limit functor, gives rise to the exact sequence of $R$-modules;
$0\lrt\underrightarrow{\lim}M_i\lrt\omega_{n_0}\lrt\omega_{n_1}\lrt\cdots,$
where $\omega_{n_i}\in\Add\omega$, meaning that $\underrightarrow{\lim}M_i\in\X_{\omega}$.
So we are done.
\end{proof}

For given two $R$-modules $M, N$, $\uHom_R(M, N)$ stands for
$\Hom_{R}(M, N)/{\mathfrak{I}(M, N)},$
where $\mathfrak{I}(M, N)$ is the $R$-submodule of $\Hom_{R}(M, N)$ consisting of all
homomorphisms factoring through a module in $\add\omega$.

As we have mentioned in the introduction, it has been proved by Chase \cite{c}
that every pure-semisimple ring is artinian. From this point of view,
the following result can be seen as a generalization of Chase's result
for the category of $\MCM$ modules. Indeed, the result below asserts that
if every balanced big $\CM$ module which belongs to $\X_{\omega}$, is fully decomposable,
 then $\underline{\Hom}_R(-, -)$ is artinian over $\MCM$ modules.

\begin{theorem}\label{prop7}Let any balanced big $\CM$ $R$-module belonging to
$\X_{\omega}$ be fully decomposable. Then $R$ is an isolated singularity.
\end{theorem}
\begin{proof}If any $\MCM$ $R$-module belongs to $\add\omega$, then $R$ will be of finite
$\CM$-type and so $R$ is known to be an isolated singularity; see \cite{a11, hl}.
Hence, in this case the desired result is obtained.
So assume that there are some $\MCM$ $R$-modules which are not in $\add\omega$.
First we show that for any $\MCM$ module $M$, $\uHom(M, M)$ is an artinian $R$-module.
To this end, clearly we only need to treat with those (indecomposable)
modules which do not belong to $\add\omega$. Suppose that $M_0$ is an arbitrary indecomposable
$\MCM$ $R$-module which does not belong to $\add\omega$.
Taking an arbitrary $R$-regular element $x$, we may have the following pull-back diagram;
{\footnotesize{\[\xymatrix@C-0.5pc@R-.8pc{&&0\ar[d]&0 \ar[d] & &&\\&&\omega^{n_1}\ar[d]\ar@{=}[r]&\omega^{n_1} \ar[d] & && \\ 0 \ar[r]  & M_0\ar@{=}[d]\ar[r]^{x} & M_1 \ar[r] \ar[d] & G' \ar[d] \ar[r]  & 0&\\
 0\ar[r]& M_0 \ar[r]^{x} & M_0 \ar[r] \ar[d] & M_0/{xM_0}\ar[r]
  \ar[r] \ar[d] & \\ &&0& 0& &&\\}\]}}
where $G'\lrt M_0/{xM_0}$ is a right minimal $\MCM$-approximation. As $\Ext_R^1(M_0, \omega^{n_1})=0$,
the left column will be split, implying that $M_1\cong M_0\oplus\omega^{n_1}$.
Evidently $M_1$ is also a $\MCM$ $R$-module which does not belonging to $\add\omega$. Thus
replacing $M_0$ by $M_1$ in the above argument, gives rise to the existence
of $R$-homomorphism $M_1\st{x}\lrt M_2$ such that $M_2\cong M_1\oplus\omega^{n_2}$.
By repeating this procedure, we obtain a chain of $R$-homomorphisms of $\MCM$ modules;
\begin{equation}\label{eq5}
M_0\st{x}\lrt M_1\st{x}\lrt M_2\st{x}\lrt\cdots,
\end{equation}
such that for any $i\geq 1$,
there is an isomorphism $M_i\cong M_0\oplus\omega^{n_i}$, for some non-negative integer $n_i$.
Applying the functor $\uHom_{R}(M_0, -)$ to \ref{eq5},
yields the following chain of $R$-modules;
\begin{equation}\label{eq6}
\uHom_{R}(M_0, M_0)\st{x}\lrt\uHom_R(M_0, M_0)\st{x}\lrt\uHom_{R}(M_0, M_0)\lrt\cdots,
\end{equation}
because $\uHom_R(M_0, M_i)\cong\uHom_R(M_0, M_0)$, for any $i\geq 1$.
According to our construction, $M_i\st{x}\lrt M_{j}$, where $j>i$,  is an $R$-monomorphism
such that its cokernel is $\MCM$ and so it can be easily seen that $\underrightarrow{\lim}M_i$ is a balanced big $\CM$ $R$-module.
In view of Lemma \ref{prop6}, $\underrightarrow{\lim}M_i\in\X_{\omega}$ and
so by the hypothesis, $\underrightarrow{\lim}M_i=\oplus_{j\in J}C_j$, where each $C_j$ is finitely generated.
Therefore, we have the following isomorphisms;
\[\begin{array}{lllll}
\underrightarrow{\lim}{\uHom_R(M_0, M_i)} & \cong \uHom_R(M_0, \underrightarrow{\lim}M_i)\\
& \cong \uHom_R(M_0, \oplus_{j\in J}C_j)\\
& \cong \oplus_{j\in J}\uHom_R(M_0, C_j).
\end{array}\]
Since $M_0$ is an indecomposable $R$-module, $\uHom_R(M_0, M_0)$
is indecomposable as an $\uEnd_R(M_0)$-module.
Now one may apply Lemma \ref{lem8} and conclude that
after some steps, all of morphisms in (\ref{eq6}) are isomorphism
or $\underrightarrow{\lim}{\uHom}_R(M_0, M_0)=0$. The former one cannot
take place. {Otherwise, we will have the isomorphism ${\uHom}_R(M_0, M_0)\cong x^t{\uHom}_R(M_0, M_0)$,
for some integer $t>0$, and so by Nakayama's lemma ${\uHom}_R(M_0, M_0)=0$, guaranteeing that
$M_0$ lies inside $\add\omega$, which is a contradiction.}
Hence the latter one will take place, meaning that there exists an integer
$t>0$ such that $x^t\uHom_R(M_0, M_0)=0$. Next suppose that
$\xx=x_1, \cdots, x_d$ is a system of parameters of $R$. As any permutation
of $\xx$ is again $R$-regular sequence, one may find an integer $n>0$
such that  $x_i^n\uHom_R(M_0, M_0)=0$, for any $1\leq i\leq d$,
that is to say, $\xx^n\uHom_R(M_0, M_0)=0$. This would imply that
$\m^u\uHom_R(M_0, M_0)=0$ for some integer $u>0$, and so $\uHom_R(M_0, M_0)$ is an artinian
$R$-module, as claimed. Next we show that the module $M_0$ is locally free on the punctured spectrum of $R$.
As we have already showed, $\uHom_R(M_0, M_0)_{\p}=0$,
for all nonmaximal prime ideals $\p$ of $R$.
Now, if $R$ is Gorenstein, i.e. $R=\omega$, then the equality
$\uHom_R(M_0, M_0)=\underline{\Hom}_R(M_0, M_0)$ gives the desired result.
In the case $R$ is not necessarily Gorenstein, it will not
belong to $\add\omega$. Thus, by repeating the above argument for $R$ instead of $M_0$,
we deduce that $\omega_{\p}=R_{\p}$, for all nonmaximal prime ideals
$\p$ of $R$, meaning that $R$ is locally Gorenstein and
consequently, ${M_0}_{\p}$ is a free $R_{\p}$-module.
{Hence any $\MCM$ $R$-module is locally free on the punctured spectrum of $R$,
and so, \cite[Lemma 3.3]{yo} yields that $R$ is an isolated singularity.}
The proof is now completed.
\end{proof}

The result below, is an immediate consequence of Theorem \ref{prop7}.

\begin{cor}\label{prop9}Let $(R, \m)$ be a complete Gorenstein local ring. If
every Gorenstein projective $R$-module is fully decomposable,
then $R$ is an isolated singularity.
\end{cor}

Let $M, N$ be two $\MCM$ $R$-modules. Recall that ${\rad}(M, N)$ is a submodule of $\Hom_R(M, N)$
consisting of those homomorphisms $\varphi:M\lrt N$ such that, when we decompose
$M=\oplus_{j}M_j$ and $N=\oplus_iN_i$ into indecomposable modules, and accordingly
decompose $\varphi=(\varphi_{ij}:M_j\lrt N_i)$, no $\varphi_{ij}$ is an isomorphism.
Moreover, ${\rad}^2(M, N)$ is  a submodule of $\Hom_R(M, N)$
consisting of those homomorphisms $\varphi:M\lrt N$ for which there is a factorization
{\footnotesize{\[\xymatrix{M\ar[dr]_{\alpha} \ar[rr]^{\varphi} && N \\ & X \ar[ur]_{\beta} & }\]}}
with $X$ is an $\MCM$ $R$-module, $\alpha\in{\rad}(M, X)$ and $\beta\in{\rad}(X, N)$. For $n>2$,
${\rad}^n(M, N)$ is defined inductively; see \cite[Definition 12.20]{lw1}.

Recall that a subcategory $\C$ of $R$-modules is
{\em of finite type,} if it has only finitely many isomorphism
classes of indecomposable modules.

The proof of the next  result is the same as in the setting
of artin algebras \cite[Lemma 3.14]{as1}, and we include it only for the sake completeness.

\begin{prop}\label{prop3}Let $R$ be an isolated singularity
and let $\C$ be a subcategory consisting of indecomposable $\MCM$ modules
that is closed under isomorphism. Let $\A$ be the subcategory consisting of
all indecomposable $\MCM$ modules that do not belong to $\C$ and assume that $\A$ is of finite type.
Then for any $X\in\add\A$ and {a} faithful system of parameters $\xx$ for $\A$, there is an $R$-homomorphism
$f:X\st{\tiny {\left[\begin{array}{ll} f_1 \\ f_2 \end{array} \right]}}\lrt M\oplus N$,
where $M\in\add\C$ and $N\in\add\A$ such that for any $L\in\add\C$, $\Hom_R(f,L)$ is surjective and
$f_2\otimes_R R/\xx^2 R=0$.
\end{prop}
\begin{proof}
We show by induction that for any $n\geq 0$ and each
$X\in\add\A$, there is a morphism
$f:X\st{\tiny {\left[\begin{array}{ll} f_1 \\ f_2 \end{array} \right]}}\lrt M\oplus N$,
where $M\in\add\C$ and $N\in\add\A$ such that for any
$L\in\add\C$, $\Hom_R(f, L)$ is surjective and $f_2\in{\rad}^n(X, N)$.
In case, $n=0$, we set $f=id_{X}$.
Now assume that $n>0$ and the result has been proved for values smaller than $n$.
By the induction hypothesis, there exists a morphism
$f':X\st{\tiny {\left[\begin{array}{ll} f'_1 \\ f'_2 \end{array} \right]}}\lrt M'\oplus N'$,
in which $f'_2\in{\rad}^{n-1}(X, N')$, $N'\in\add\A, M'\in\add\C$
and for any $K\in\add\C$, the morphism $\Hom_R(f', K)$ is surjective.
Since the category of $\MCM$ modules has left almost split morphisms,
there is an $R$-homomorphism
$g:N'\st{\tiny {\left[\begin{array}{ll} g_1 \\ g_2 \end{array} \right]}}\lrt Z=Z_1\oplus Z_2$,
with $Z_1\in\add\A$, $Z_2\in\add\C$ and $g_1\in{\rad}(N', Z_1)$ such that $\im\Hom(g, L)=\Hom_R(N', L)$
for any $L\in\add\C$.
Assuming $h$ as the following composition morphism $$h:X\st{\tiny {\left[\begin{array}{ll} f'_1 \\ f'_2 \end{array} \right]}}\lrt
M'\oplus N'\st{id_{M'}\oplus g}\lrt M'\oplus Z$$
where $g_1f'_2\in {\rad}^n(X, Z_1)$, we have that $\Hom(h, K)$ is surjective for any $K\in\add\C$.
On the other hand, suppose that $\F:=\{X_1, \cdots, X_t\}$ is the set of all pairwise
non-isomorphic indecomposable objects of $\A$.
{Since $\xx$ is} a faithful system of parameters for $\A$,
$\F/\xx^2\F=\{X_1/\xx^2X_1 \cdots, X_t/\xx^2X_t\}$
is a set of indecomposable modules of finite length.
By virtue of Corollary to \cite[Lemma 12]{hs}, there is a non-negative integer $n$ such
that ${\rad}^n(X_i/{\xx^2 X_i}, X_j/{\xx^2 X_j})=0$ for all $X_i, X_j\in\F$.
Consequently, $\bar{f}=f\otimes_R R/\xx^2 R\in{\rad}^n(X_i/\xx^2X_i,X_j/\xx^2 X_j)$ and
so $\bar{f}=0$,
which gives the desired result.
\end{proof}

\begin{theorem}\label{fcmt}Let $R$ be an isolated singularity which
is not of finite $\CM$-type. Then there is an infinite
set of pairwise non-isomorphic indecomposable $\MCM$ modules $\{M_i\}_{i\in\N}$ and
non-zero $R$-homomorphisms $f_i:M_i\lrt k$ such that any composition map
$M_j\lrt M_i\st{f_i}\lrt k$ is zero, for all $j>i$.
\end{theorem}
\begin{proof}
In order to obtain the desired result,
we will first construct a pairwise disjoint
infinite family of finite type subcategories of indecomposable $\MCM$ modules
$\A_1, \A_2, \A_3, \cdots,$ inductively as follows:\\
$(i)$ We set $\mathcal{A}_1$ to be the class of all projective $R$-modules that are isomorphic to $R$.
$(ii)$ Suppose $j>1$ is an integer and assume that we have already
constructed $\A_1, \cdots, \A_{j-1}$. Letting $\C_j=\ind(\MCM)-\bigcup_{i=1}^{j-1}\A_i$
{and $\xx$ a faithful system of parameters for $\bigcup_{i=1}^{j-1}\A_i$,}
by Proposition \ref{prop3}, there is an $R$-homomorphism
$f:R\st{\tiny {\left[\begin{array}{ll} f_1 \\ f_2 \end{array} \right]}}\lrt K_j\oplus N$,
where $K_j\in\add\C_j$ and $N\in\add(\bigcup_{i=1}^{j-1}\A_i)$
such that for any $L\in\add\C_j$, $\Hom_R(f,L)$ is surjective and
$f_2\otimes_R R/\xx^2 R=0$. 
By the Krull-Remak-Schmidt theorem, $K_j=\oplus_{i=1}^t X_i$, where each $X_i$
is an indecomposable finitely generated submodule of $K_j$. We put $\mathcal{A}_j$ to be the class of all $\MCM$ modules
that are isomorphic to one of $X_1, \cdots, X_t$.\\
So we have constructed a pairwise disjoint infinite family of finite type
subcategories of indecomposable $\MCM$ modules $\A_1, \A_2, \cdots$.\\
Let us divide the remainder of the proof into three steps:\\
{\sc step 1}: We show that, for any $j$, $\mathcal{A}_j$ is a generator for $\mathcal{C}_s$, for any
$s\geq j$, namely, for each $L\in\mathcal{C}_s$, there exists an $R$-epimorphism $Y\lrt L$, where $Y\in\add\A_j$.
To see this, take an arbitrary object $L\in\mathcal{C}_s$ and consider an epimorphism $\alpha:R^n\lrt L$.
By part ($ii$), there exists an $R$-homomorphism $f:R\st{\tiny {\left[\begin{array}{ll} f_1 \\ f_2 \end{array} \right]}}\lrt K_j\oplus N$ such that $K_j\in\add\mathcal{C}_j$ and $N\in\add (\bigcup_{i=1}^{j-1}\mathcal{A}_i)$.
Since for any $s\geq j$, $\mathcal{C}_s\subseteq\mathcal{C}_j$, $L\in\mathcal{C}_j$ and by construction of
$f$ in part $(ii)$, $\Hom_R(f, L)$ is surjective. Thus there is an $R$-homomorphism $(\psi_1, \psi_2):K_j^n\oplus N^n\lrt L$
such that the diagram
{\footnotesize{\[\xymatrix{R^n\ar[dr]_{\alpha} \ar[rr]^{\tiny {\left[\begin{array}{ll} f_1^n \\ f_2^n \end{array} \right]}} &&\ar[dl]^{(\psi_1, \psi_2)} K_j^n\oplus N^n \\ & L& }\]}}
is commutative. {Take a faithful system of parameters $\xx$ for $\bigcup_{i=1}^{j-1}\mathcal{A}_i$,
which has been used in the construction of $\A_j$ in part $(ii)$.}
Now, by applying the functor $-\otimes_R R/\xx^2 R$ and using the fact that $f_2\otimes_R R/\xx^2 R=0$,
we infer that $\bar{\psi_1}$ : $K_j^n/\xx^2 K_j^n\lrt L/\xx^2L$ is an epimorphism and so
by Nakayama's lemma $\psi_1$, will be an epimorphism, as well. We set $Y:= K_j^n$.\\
{\sc step 2}: Next we show that for any $j$, there is an object $X_j$ of $\A_j$  such that
there is not any epimorphism $\oplus_{i>j}M_i^{(s_i)}\lrt X_j$
where $M_i^,$s are objects of $\A_i^,$s and $s_i$ is a set for each $i$.
Suppose that for some $j$, this is not the case.
Let $\A_j = \{X_1, X_2, \cdots, X_t\}$, up to isomorphism. As $\oplus_{i=1}^tX_i$ is finitely generated,
we may assume that there is an epimorphism $\varphi: \oplus_{i=j+1}^nM_{i}^{n_i}\lrt \oplus_{i=1}^tX_i$ for some integer $n>0$,
where $M_i^,$s belong to $\A_i^,$s and $n_i>0$ is an integer for each $i$. Since $R$ is an isolated singularity, there exists
a faithful system of parameters $\yy$ for $\{M_{j+1}, \cdots, M_{n}, X_1, \cdots, X_t\}$, by Proposition \ref{cor1}.
Considering  the  epimorphism $\bar{\varphi}:\oplus_{i=j+1}^n(M_{i}^{n_i}/{\yy^2M_{i}^{n_i}})\lrt (\oplus_{i=1}^tX_i)/\yy^2(\oplus_{i=1}^tX_i)$,
Main property (a) of \ref{s1} yields that $${\text {min}}
\{ \mu^*(X_i/{\yy^2X_i})| 1\leq i\leq t \} \geq {\text {min}}
\{\mu^*(M_{i}/{\yy^2M_{i}}) | j+1\leq i \leq n \}.$$
It should be observed that if the equality takes place, then
by part $(b)$ of Main property of \ref{s1}, $\bar{\varphi}$ is a split epimorphism. On the other hand by \cite[Corollary 15.11]{lw1},
quotient modules, for any $1\leq i\leq t$, $X_i/\yy^2 X_i$ and for any $j+1\leq i\leq n$, $M_i/\yy^2 M_i$
are indecomposable. Thus by the Krull-Remak-Schmidt theorem,
for some $1\leq i\leq t$, $X_i/\yy^2 X_i$ is isomorphism with a direct summand of $\oplus_{i=j+1}^n M_j/\yy^2 M_j$,
{and so, $X_i/\yy^2 X_i\cong M_s/\yy^2 M_s$ for some integer $j+1\leq s\leq n$.
Hence, applying \cite[Lemma 3.3.2]{bh} gives rise to the isomorphism $X_i\cong M_s$, which contradicts with our construction of $\A_i^,$s.}
Assuming $\mu^*(M_{s}/{\yy^2M_{s}})={\text {min}}\{\mu^*(M_{i}/{\yy^2M_{i}})\mid j+1\leq i\leq n\}$,
by the step 1, there is an epimorphism
{$\oplus_{i=1}^tX_i^{m_i}\lrt M_{s}$, where $m_i>0$ is an integer for any $i$ and so)}
$(\oplus_{i=1}^tX_i^{m_i}/{\yy^2(\oplus_{i=1}^tX_i^{m_i}}))\lrt M_{s}/{\yy^2M_{s}}$ will be an
epimorphism, meaning that
$\mu^*(M_{s}/{\yy^2M_{s}})>\mu^*(X_j/{\yy^2X_j})$ for some $1\leq j\leq t$, however, this is impossible.\\
{\sc step 3}:
{As it has been observed in the previous step, for any $i$, there exists
an indecomposable object $M_i\in\A_i$ such that there is not any epimorphism
$\varphi:\oplus_{j=i+1}^nM_j^{n_j}\lrt M_i$, where $M_j^,$s are indecomposable objects of $\A_j^,$s.}
 Suppose that $\F$ is the class consisting of all these $M_i^,$s.
We should stress that, also there may exist more than one module $M_i\in\A_i$, with
the mentioned property; however, we put only one of them in $\F$.
Since $\mathcal{A}_i^,$s are pairwise disjoint infinite family,
$\F$ will be an infinite set of indecomposable pairwise non-isomorphic $\MCM$ modules.
The same argument given in the proof of step 3 of Theorem \ref{the2},
ensures the existence of non-zero $R$-homomorphisms $f_i:M_i\lrt k$ such that
any composition map $M_j\lrt M_i\st{f_i}\lrt k$ with $j>i$, is zero.
The proof then is finished.
\end{proof}

Now we are in a position to state and prove the main result of this section,
which is presented as Theorem \ref{fcmti} in the introduction.

\begin{theorem}\label{th6}The following conditions are equivalent:
\begin{enumerate}
\item $R$ is of finite $\CM$-type.
\item
The subcategory of $R$-modules consisting of all balanced big $\CM$ modules $M$ having
an $\m$-primary cohomological annihilator coincides with $\fd$.
\item
$R$ is an isolated singularity and the class $\fd$ is closed under kernels of epimorphisms.
\item
$R$ is an isolated singularity and $\fd$ is closed under extensions and direct summands.
\end{enumerate}
\end{theorem}
\begin{proof}
(1) $\Rightarrow$ (2): 
Theorem \ref{th7} gives the desired result.

(2) $\Rightarrow$ (3): Take an arbitrary $\MCM$ $R$-module $M$. As $M$ lies in $\fd$, by the hypothesis,
$M$ has an $\m$-primary cohomological annihilator. So, $R$ is an isolated singularity.
Next consider a short exact sequence of $R$-modules;
$0\lrt M'\lrt M\lrt M'' \lrt 0$, where $M, M''$ belong to $\fd$. By the hypothesis,
$M, M''$ are balanced big $\CM$ modules with $\m$-primary cohomological annihilators,
implying that $M'$ is a weak balanced big $\CM$ $R$-module, and by 
{Remark \ref{re1},}
we get that $M'$ has an $\m$-primary cohomological annihilator.
Hence invoking Lemma \ref{lem1}, yields that either $M'$ is zero or it is balanced big
$\CM$. Consequently, $M'$ is in $\fd$.

(3) $\Rightarrow$ (1): Assume on the contrary that $R$ is not of finite $\CM$-type.
So in view of Theorem \ref{fcmt}, there is an infinite set of pairwise non-isomorphic
indecomposable $\MCM$ $R$-modules $\{M_i\}_{i\in I}$ and non-zero $R$-homomorphisms $f_i:M_i\lrt k$ such that
any composition map $M_j\lrt M_i\st{f_i}\lrt k$ with $j>i$, is zero. Here $I$ is a subset of $\N$.
As $G\st{\alpha}\lrt k$ is a right minimal $\MCM$-approximation, there is a non-zero
homomorphism $g_i:M_i\lrt G$, for any $i$ such that $\alpha g_i=f_i$.
Setting $g=(g_i)_{i\in I}:\oplus_{i\in I} M_i\lrt G$, we have a short exact sequence of $R$-modules;
$0\lrt K\st{\theta}\lrt (\oplus_{i\in I}M_i)\oplus R^n\st{[g~~~\beta]}\lrt G\lrt 0$.
Hence, the hypothesis $\fd$ being closed under kernels of epimorphisms, yields that
$K$ is in $\fd$.  Now Proposition \ref{prop10} forces $I$ to be a finite set,
which is a contradiction. Therefore, $R$ is of finite $\CM$-type.

(1) $\Rightarrow$ (4): According to Corollary 2 of \cite{hl}, $R$ is an isolated singularity.
Next, consider a short exact sequence of $R$-modules; $0\lrt M'\lrt M\lrt M''\lrt 0$,
in which $M',M''\in\fd$. By the hypothesis, $M'$ and $M''$ are balanced big $\CM$ modules
with $\m$-primary cohomological annihilators. So it is fairly easy to see that
$M$ is balanced big $\CM$ with an $\m$-primary cohomological annihilator.
Now Theorem \ref{th7} would imply that $M\in\fd$.
This means that $\fd$ is closed under extensions.
Moreover, Theorem \ref{th2} indicates that $\fd$ is closed under direct summand.

(4) $\Rightarrow$ (3).
Take a short exact sequence of $R$-modules; $0\lrt M'\lrt M\lrt M''\lrt 0$,
in which $M, M''\in\fd$. We would like to show that $M'\in\fd$.
{By the hypothesis, $M''=\oplus_{i\in I}X_i$, where each $X_i$ is finitely generated.
For any $i\in I$, take a short exact sequence of finitely generated $R$-modules,
$0\lrt L_i\lrt P_i\lrt X_i\lrt 0$, in which $P_i\lrt X_i$ is a projective cover.
In particular, one may have the short exact sequence of $R$-modules,
$0\lrt L\lrt P\lrt M''\lrt 0$, where $P=\oplus_{i\in I}P_i$ and $L=\oplus_{i\in I}L_i$.}
Considering the following commutative diagram with exact rows;

{\footnotesize{$$\begin{CD}
0 @>>> L  @>>> P @>>> M'' @>>> 0\\
& & @Vu VV @V VV @V id VV & &\\ 0 @>>> M' @>>> M @>>> M''
@>>> 0,\end{CD}$$}}
we obtain the exact sequence $0\lrt L\lrt P\oplus M'\lrt M\lrt 0$.
Since $L, M\in\fd$, by our assumption, $P\oplus M'\in\fd$.
Consequently, $M'\in\fd$, because by the hypothesis $\fd$ is closed
under direct summand. So the proof is completed.
\end{proof}

Here we recover a notable result of Beligiannis \cite[Theorem 4.20]{be}.

\begin{theorem}\label{th5} Let $R$ be a Gorenstein complete local ring. Then
the following conditions are equivalent:
\begin{enumerate}
\item $R$ is of finite $\CM$-type.
\item
Every Gorenstein projective $R$-module is fully decomposable.
\item
The subcategory of Gorenstein projective $R$-modules with
$\m$-primary cohomological annihilators coincides with $\fd$.
\item
 The category of all indecomposable finitely generated Gorenstein projective $R$-modules is  of bounded $\uh$-length.
\end{enumerate}
\end{theorem}
\begin{proof}
(1) $\Rightarrow$ (2): This is  Theorem \ref{th8}.

(2) $\Rightarrow$ (1):
In view of Corollary \ref{prop9}, $R$ is an isolated singularity.
Now by applying the proof of the implication ($3 \Rightarrow 1$) of Theorem \ref{th6}
and using the fact that the category of Gorenstein projective modules is closed
under kernels of epimorphisms, we deduce that $R$ is of finite $\CM$-type.

(3) $\Rightarrow$ (1): By the assumption, every $\MCM$ $R$-module has an $\m$-primary cohomological
annihilator, implying that $R$ is an isolated singularity.
Moreover, it follows from the hypothesis that $\fd$ is closed under kernels of
epimorphisms. Now the implication ($3 \Rightarrow 1$) of Theorem \ref{th6} yields
the required result.

(1) $\Rightarrow$ (3): This follows from the implication (1) $\Rightarrow$ (2).

(4) $\Leftrightarrow$ (1): The implication (4) $\Rightarrow$ (1) follows from Corollary \ref{cor2}, whereas the reverse
implication holds trivially.
\end{proof}

\section{Representation properties of $\CM$ modules over artin algebras}

Motivated by (commutative) complete Cohen-Macaulay local rings,
Auslander and Reiten in \cite{AR1, AR} have introduced and studied Cohen-Macaulay artin algebras.
Recall that an artin algebra $\Lambda$ is said to be {\em Cohen-Macaulay}
if there exists a pair of adjoint functors $(G, F)$ between $\md\Lambda$ and $\md\Lambda$,
inducing mutually inverse equivalences;
\[\xymatrix{\mathcal{I}^{\infty}(\Lambda) \ar@<1ex>[rr]^{F}&& \mathcal{P}^{\infty}(\Lambda)\ar@<1ex>[ll]^{G},}\]
where $\mathcal{P}^{\infty}(\Lambda)$ (resp.  $\mathcal{I}^{\infty}(\Lambda))$ denotes the category of
all finitely generated modules of finite projective (resp. injective) dimension.

As we have noted in the introduction, it is well-known that if $\Lambda$ is a Cohen-Macaulay
artin algebra, then there is a finitely
generated $\Lambda$-bimodule $\omega$ such that the functors  $F$ and $G$ are presented by
$\Hom_{\Lambda}(\omega,-)$ and $\omega\otimes_{\Lambda}-$, respectively; see \cite{AR1}.
In this case, $\omega$ is called a {\em dualizing module} for $\Lambda$.

\begin{remark}There is a tight connection between dualizing modules and strong cotilting
modules over artin algebras. Precisely, a $\Lambda$-bimodule $\omega$ is dualzing if and
only if $\omega$ is strong cotilting viewed both as left and right modules
{and the natural map $\Lambda\lrt\End(_{\Lambda}\omega)$
is an isomorphism}; see \cite[Proposition 3.1]{AR1}.
This connection gives an interesting interplay between cotilting theory
for artin algebras and module theory for commutative Cohen-Macaulay rings.
A selforthogonal $\Lambda$-module $\omega$ is cotilting if $\id_{\Lambda}\omega<\infty$
and all injective $\Lambda$-modules are in $\widehat{\add\omega}$,
and it is said to be strong cotilting if, moreover, the equality
$\mathcal{I}^{\infty}(\Lambda)=\widehat{\mathsf{add}\omega}$ holds.
{Recall that $\omega$ is said to be {\em selforthogonal}, provided that $\Ext_{\Lambda}^i(\omega,\omega)=0$
for all $i>0$.

We emphasize that the results of this section
remain true even if $\omega$ is assumed to be a cotilting $\Lambda$-module
and the natural map $\Lambda\lrt\End(_{\Lambda}\omega)$
is an isomorphism.
We are indebted to Professor Osamu Iyama for pointing us this fact.}
\end{remark}
\begin{s}
Throughout this section, $\Lambda$ is always a Cohen-Macaulay artin algebra
and $\omega$ is a dualizing $\Lambda$-bimodule.
We say that a $\Lambda$-module $M$ is {\it $\omega$-Gorenstein projective}, if it admits
a right resolution by modules in $\Add\omega$, that is, an exact sequence of $\Lambda$-modules;
$$0\lrt M\lrt w_0\st{d_0}\lrt w_1\st{d_1}\lrt\cdots \st{d_{i-1}}\lrt w_i\st{d_i}\lrt\cdots,$$ with $w_i\in{\Add\omega}$.
So finitely generated $\omega$-Gorenstein projective modules are Cohen-Macaulay in the sense
of Auslander and Reiten \cite{AR} and we also call them Cohen-Macaulay modules ($\CM$ modules).
{It should be noted that since $\omega$ is a selforthogonal $\Lambda$-module of finite injective
dimension, $\Ext_{\Lambda}^{i>0}(W, W')=0$ for all modules $W,W'\in\Add\omega$.
{Indeed, this follows from the isomorphisms $\Ext_{\Lambda}^i(\oplus_{j\in J}\omega, W')\cong\prod_{j\in J}\Ext_{\Lambda}^i(\omega, W')$ and $\Ext_{\Lambda}^i(\omega, \oplus_{j\in J'}\omega)\cong\oplus_{j\in J'}\Ext_{\Lambda}^i(\omega, \omega)$. One should observe that,
as $\omega$  admits a projective resolution of finitely generated projective modules, \cite[Exercise 2(a), page 16]{ej}
ensures the validity of the latter isomorphism.} So it is easily seen that
our notion of $\omega$-Gorensteiness coincides with the
one given by Holm and J${\o}$rgensen in \cite{hj}.}
We say that an $\omega$-Gorenstein projective $\Lambda$-module $M$ is {\it fully decomposable} (resp. {\it of finite
$\CM$-type}) if  it is a direct sum of arbitrarily many
copies (resp. of a finite number up to isomorphisms) of indecomposable $\CM$ modules.

Moreover, by {\it $\CM$-support} of an $\omega$-Gorenstein projective module $M$,
denoted by $\CM$-$\Supp_{\Lambda}(M)$, we mean the set of all indecomposable $\CM$
$\Lambda$-modules $N$ such that $\Hom_{\Lambda}(N,M)\neq 0$.
\end{s}

Our aim in this section is to examine results in the previous sections
in the context of Cohen-Macaulay artin algebras.
It is proved that any $\omega$-Gorenstein projective $\Lambda$-module with bounded length on $\CM$-support
must be fully decomposable; see Theorem \ref{the9}. In particular, it will be observed
in Theorem \ref{th4} that if an
$\omega$-Gorenstein projective module  $M$ is not of finite $\CM$-type,
then there are indecomposable $\CM$ $\Lambda$-modules of arbitrarily large
(finite) length, guaranteeing the validity of the
first Brauer-Thrall conjecture for the category of  Cohen-Macaulay modules over
Cohen-Macaulay artin algebras. Moreover, our results  extend
a result of Chen \cite[Main Theorem]{che} for Cohen-Macaulay artin
algebra, that is, we specify Cohen-Macaulay
artin algebras of finite $\CM$-type in terms of the decomposition properties of
$\omega$-Gorenstein projective modules.

Let $\Lambda\ltimes\omega$ denote the trivial extension of $\Lambda$ by $\omega$.
Then according to ring homomorphisms; $\Lambda\lrt\Lambda\ltimes\omega\lrt\Lambda$,
any $\Lambda$-module can be viewed as a $\Lambda\ltimes\omega$-module and vise versa,
and in this section we shall do so freely.

Assume that $\F$ is a class of $\Lambda$-modules and $M$ a $\Lambda$-module.
A homomorphism $f:M\lrt F$, where $F\in\F$, is said to be an $\F$-{\it preenvelope} of $M$,
provided that for every homomorphism $g:M\lrt F'$, where $F'\in\F$,
there exists a homomorphism $h:F\lrt F'$ such that $hf=g$.

\begin{prop}\label{lem4}Every $\omega$-Gorenstein projective $\Lambda$-module is a direct limit of $\CM$ modules.
\end{prop}
\begin{proof}Take an arbitrary $\omega$-Gorenstein projective $\Lambda$-module $M$.
Because of \cite[Proposition 2.1]{le}, it suffices to show that
any $\Lambda$-homomorphism $f:N\lrt M$, where $N$ is finitely generated,
factors through a $\CM$ $\Lambda$-module, say $C$. {Assume that $\id_{\Lambda}\omega=n$}. In view of \cite[Proposition 2.13]{hj},
$M$ is Gorenstein projective over $\Lambda\ltimes\omega$, so one may take the following
exact sequence of $\Lambda\ltimes\omega$-modules;
$$0\lrt M\lrt Q^0\lrt\cdots\lrt Q^{n-1}\lrt L\lrt 0,$$ in which for any
$i$, $Q^i$ is projective  and $L$ is Gorenstein projective.
Another use of \cite[Proposition 2.13]{hj} yields that as a $\Lambda$-module,
$L$ is $\omega$-Gorenstein projective.
Since $N$ is a finitely generated $\Lambda$-module, evidently it is
finitely generated over $\Lambda\ltimes\omega$.
Consider the following sequence of finitely generated $\Lambda\ltimes\omega$-modules;
$$N\st{d^0}\lrt P^0\st{d^1}\lrt\cdots\lrt P^{n-1}\lrt K\lrt 0,$$
where $N\lrt P^0$ and $\cok(d^i)\lrt P^{i+1}$, for any $i$,
are projective preenvelopes. {It should be noted that these preenvelopes exist
because of \cite[page 247]{ej}.} According to \cite[Theorem 4.32]{fgr}, $\Lambda\ltimes\omega$
is a Gorenstein algebra with injective dimension $n$, where $n=\id_{\Lambda}\omega.$
Hence by using \cite[Theorem 10.2.14]{ej}, we have the following exact sequence of finitely generated
$\Lambda\ltimes\omega$-modules; $$0\lrt C\lrt P_{n-1}\lrt\cdots\lrt P_0
\lrt K\lrt 0,$$ such that each $P_i$ is projective and $C$
is Gorenstein projective.
Consequently, one obtains the following commutative
diagram of $\Lambda\ltimes\omega$ (and also $\Lambda$)-modules;
 which is similar to diagram appeared in the proof of \cite[Lemma 10.3.6]{ej};
{\tiny{\[\xymatrix@C-0.6pc@R-1.4pc{ && N \ar[rr] \ar[ddd]_{f}
\ar[ddr]^{h_n}  && P^0 \ar[rr] \ar[ddd]_{f_0} \ar[ddr]^{h_{n-1}} &&  \cdots \ar[rr] && K
\ar[rr] \ar[ddr]^{id} \ar[ddd]_{f_n} && 0 \\ \\
& 0 \ar[rr] && C \ar[rr] \ar[dl]^{g_n} &&  P_{n-1} \ar[rr] \ar[dl]^{g_{n-1}} &&
\cdots \ar[rr] && K \ar[rr] \ar[dl]^{f_n} && 0 \\
0 \ar[rr] && M \ar[rr] && Q^0 \ar[rr] && \cdots \ar[rr] && L \ar[rr]&& 0.
}\]}}
One should note that the morphisms $h_i^,$s are lifted from $id:K\lrt K$,
whereas, the existence of $f_i^,$s follows from the construction of upper row.
Finally, the morphisms $g_i^,$s exist, because they are lifted from
$f_n$. Now chasing diagram enables us to deduce that $f$ factors through $C\oplus P^0$
which is $\CM$ over $\Lambda$, thanks to \cite[Proposition 2.13]{hj}.
Clearly $f$ factors from this module as a $\Lambda$-homomorphism.
So the proof is finished.
\end{proof}

We need the following result for later use.
\begin{lemma}\label{lem3}Let $M$ be a non-zero $\omega$-Gorenstein projective $\Lambda$-module.
If $\CM$-$\Supp_{\Lambda}(M)$ is of bounded length, then $M$ has an indecomposable
$\CM$ direct summand.
\end{lemma}
\begin{proof}
According to Proposition \ref{lem4}, there is a direct system
of $\CM$ $\Lambda$-modules $\{M_i, \varphi^i_j\}_{i, j\in I}$ such that
$M=\underrightarrow{\lim}M_i$.
As $M$ is non-zero, we can take an index $j\in I$ and an indecomposable $\CM$ direct summand $X_j$ of $M_j$
such that the morphism ${\varphi_j}_{\mid_{X_j}}:X_j\lrt M$ is non-zero, 
{where $\varphi_j:M_j\lrt M$ is the natural morphism such that for any $i\leq j$, $\varphi_i=\varphi_j\varphi^i_j$}. Let $k_1\in I$ be an index with
$k_1> j$. So we have an indecomposable $\CM$ direct summand $X_{k_1}$ of $M_{k_1}$ such that
$$X_j\st{\varphi_{k_1}^j|_{X_j}}\lrt M_{k_1}\st{\pi}\lrt X_{k_1}\st{{\varphi_{k_1}}{|_{X_{k_1}}}}\lrt M$$
is non-zero, where $\pi:M_{k_1}\lrt X_{k_1}$ is the canonical projection. We denote the composition map
$X_j\st{\varphi_{k_1}^j{|_{X_j}}}\lrt M_{k_1}\st{\pi}\lrt X_{k_1}$ by $\psi_{k_1}^j$. One can use the
induction argument to obtain a chain of morphisms of indecomposable $\CM$ $\Lambda$-modules
$$X_j\st{\psi_{k_1}^j}\lrt X_{k_1}\st{\psi_{k_2}^{k_1}}\lrt X_{k_2}\st{\psi_{k_3}^{k_2}}\lrt X_{k_3}
\lrt\cdots,$$
such that {any composite of finite
number of morphisms has non-zero image in $M$.}
Since all $X_i^,$s belong to $\CM$-$\Supp_{\Lambda}(M)$,
they are of bounded length. Hence, Harada-Sai Lemma \cite[Lemma 11]{hs}, guarantees the existence
of an index $k_t\in I$ such that for each $k_s> k_t$, the induced morphism $\psi^{k_t}_{k_s}$
needs to be an isomorphism. This implies that, for any $k_s> k_t$ the morphism
$\varphi_{k_s}^{k_t}|_{X_{k_t}}:X_{k_t}\lrt M_{k_s}$ is a split monomorphism. This, in turn,
would imply that ${\varphi_{k_t}}_{\mid_{X_{k_t}}}:X_{k_t}\lrt M$ is a pure monomorphism. As $X_{k_t}$ is a finitely generated
module over the artinian ring $\Lambda$, it will be pure injective, inforcing ${\varphi_{k_t}}_{\mid_{X_{k_t}}}$ to be  a
split monomorphism. Hence, $M$ has an indecomposable $\CM$ direct summand $X_{k_t}$.
So the proof is finished.
\end{proof}


The next result indicates that, for a given $\omega$-Gorenstein projective module
$M$, the boundedness of its $\CM$-support forces $M$ to be fully decomposable.

\begin{theorem}\label{the9}
Let $M$ be an $\omega$-Gorenstein projective $\Lambda$-module. If $\CM$-$\Supp_{\Lambda}(M)$
is of bounded length, then $M$ is fully decomposable.
\end{theorem}
\begin{proof}
According to Lemma \ref{lem3}, $M$ has an indecomposable $\CM$ direct summand $X$.
{Put $\Sigma$ to be the set of all fully decomposable pure submodules of $M$}  
 For any two objects $N, L\in\Sigma$, we write
$N\leq L$ if and only if 
the following diagram
is commutative;
{\footnotesize{\[\xymatrix@C-0.5pc@R-.8pc{N\ar[dr]_{i_N} \ar[rr]^{i_{NL}} && L\ar[dl]^{i_L} \\ & M  & }\]}}
where $i_N, i_L, i_{NL} $ are pure monomorphism.
 \new{and $i_{NL}$ is the inclusion map.}
{Assume that $Y=\oplus X_i$ is a pure submodule of} $M,$ where each
$X_i$ is an indecomposable $\CM$ direct summand of $M$, and $Y$ is maximal with respect
to this property. Take the pure exact sequence of $\Lambda$-modules;
$$\eta:0\lrt Y\st{i_Y}\lrt M\lrt K\lrt 0.$$
{Let $f:N\lrt K$ be a non-zero $\Lambda$-homomorphism, where $N$ is finitely generated.
As $\eta$ is a pure exact sequence, $f$ will factor through $M$.  In view of Proposition \ref{lem4},
$M=\underrightarrow{\lim}M_i$, where each $M_i$ is a $\CM$ module. Consequently, for some index $i$,
the morphism $f$ factors through $M_i$, and so, \cite[Proposition 2.1]{le}, enables us to infer that $K=\underrightarrow{\lim}K_i,$
where each $K_i$ is a $\CM$ $\Lambda$-module. On the other hand,}
it is evident that any element of $\CM$-$\Supp_{\Lambda}(K)$ belongs to
$\CM$-$\Supp_{\Lambda}(M)$, and so $\CM$-$\Supp_{\Lambda}(K)$ will be of bounded length.
{Therefore, by virtue of Lemma \ref{lem3}, $K$ has an indecomposable $\CM$ direct summand $X$.
Thus $Y\oplus X$ is a pure submodule $M$, containing $Y$ properly. However, this contradicts the maximality of $Y$.
Hence, $K=0$, and then we get the isomorphism $Y\cong M$. So the proof is completed.}
\end{proof}

\begin{theorem}\label{th4}Let $M$ be an $\omega$-Gorenstein projective $\Lambda$-module
which is not of finite $\CM$-type. Then there are indecomposable $\CM$ $\Lambda$-modules of arbitrarily
large finite length.
\end{theorem}
\begin{proof}Assume for the contradiction that the class of all indecomposable
$\CM$ $\Lambda$-modules is of bounded length.
So by Theorem \ref{the9}, we deduce that $M$ is fully decomposable.
Suppose that  $M=\oplus_{i\in I}M_i^{(t_i)}$, in which
for any $i$, $M_i$ is an indecomposable $\CM$ $\Lambda$-module. Put $\F=\{M_i\mid i\in I\}.$
By our assumption, $\F$ is of bounded length. By property 2 of \ref{s1}, there are only finitely many
Gabriel-Roiter comeasures for $\F$. Thus it is not a restriction if we
additionally  assume that all modules in $\F$ have a fixed Gabriel-Roiter comeasure.
Suppose that $\{S_1, \cdots, S_n\}$ is the complete list of non-isomorphic simple
$\Lambda$-modules. Putting $S=\oplus_{j=1}^nS_j$, analogous to the proof of Theorem \ref{the2}
(steps 2 and 3), for each $i$,
there is a $\Lambda$-homomorphism $f_i:M_i\lrt S$ such that for any $j\in I$
with $i\neq j$, any composition map $M_i\lrt M_j\st{f_j}\lrt S$ is zero.
In view of \cite[Proposition 1.4]{AR1}, there exists a right $\CM$-approximation
$\alpha' :G'\lrt S$, and so for any $i\in I$,
one may find a $\Lambda$-homomorphism $g_i:M_i\lrt G'$ such that $\alpha' g_i=f_i$.
Set $g=(g_i)_{i\in I}:\oplus_{i\in I}M_i\lrt G'$.
Consider the exact sequence of $\Lambda$-modules;
$0\lrt K\st{\theta}\lrt (\oplus_{i\in I}M_i)\oplus \Lambda^n\st{[g~~~\beta]}\lrt G'\lrt 0$.
Evidently, $K$ is $\omega$-Gorenstein projective and hence any direct summand of $K$ is again
$\omega$-Gorenstein projective. Consequently, by Theorem \ref{the9},
$K=\oplus_{i\in J}K_i$, where for any $i$, $K_i$ is an indecomposable
$\omega$-Gorenstein projective $\Lambda$-module.
Now a similar result to Proposition \ref{prop10} leads us to infer that
$I$ is a finite set, meaning that $M$ is of finite $\CM$-type.
\end{proof}

The result below, which is an immediate consequence of Theorem \ref{th4}, should be seen
as the first Brauer-Thrall theorem for $\CM$ $\Lambda$-modules.

\begin{cor}\label{cor10}Let the category of all indecomposable $\CM$ $\Lambda$-modules
be of bounded length. Then $\Lambda$ is of finite $\CM$-type.
\end{cor}

\begin{s}\label{s5}According to \cite{AR1}, the category of $\CM$ $\Lambda$-modules
admits almost split sequences. Moreover, for a given object $M\in\md\Lambda$,
by \cite[Proposition 1.4]{AR1} there is a $\CM$-approximation $X\lrt M$. Hence
one may deduce that the category of $\CM$ modules has left almost split morphisms.
Assume that $\A$ is a finite type subcategory of
$\CM$ $\Lambda$-modules. Then the same argument given in the proof of Proposition
\ref{prop3} (see also \cite[Proposition 3.13]{as1}) indicates that
for any $X\in\add\A$, there is a $\Lambda$-homomorphism $f:X\lrt M$ with $M\in\add\C=\CM-\A$
such that for any $L\in\add\C$, $\Hom(f, L)$ is surjective, that is to say,
$f:X\lrt M$ is a $\C$-preenvelope.
\end{s}

\begin{theorem}\label{th3}
Every $\omega$-Gorenstein projective module is fully decomposable if
and only if $\Lambda$ is of finite $\CM$-type.
\end{theorem}
\begin{proof}  The `if' part is Theorem \ref{the9}.
For the `only if' part, assume that $\Lambda$ is not of finite $\CM$-type.
Analogously to the proof of Theorem \ref{fcmt},
we obtain a pairwise disjoint infinite family of finite type subcategories of indecomposable $\CM$ modules
$\A_1, \A_2, \A_3, \cdots$ as follows:\\
$(i)$ Assume that $\A_1$ is the class of all indecomposable $\Lambda$-modules
that are isomorphic to indecomposable projective $\Lambda$-modules.\\
$(ii)$ Suppose that for any $j>1$, we have already constructed $\mathcal{A}_1, \cdots, \mathcal{A}_{j-1}$.
Set $\C_j=\ind\CM-(\bigcup_{i=1}^{j-1}\A_i)$.
For a given $Q\in\Add\A_1$,
take a $\C_j$-preeenvelope $f:Q\lrt K_j$, which exists by \ref{s5}.
By the Krull-Remark-Schmidt theorem, $K_j=\oplus_{i=1}^tX_i$, where each $X_i$
is a finitely generated indecomposable submodule of $K_j$. Now put $\A_j$ to be the class
of all $\CM$ modules that are isomorphic to one of $X_1, \cdots, X_t$.\\
Let us divide the remainder of the proof into three steps:\\
{\sc step 1}: We show that, for any $j$, $\mathcal{A}_j$ is a generator for $\mathcal{C}_s$,
for any $s\geq j$. To see this,
take an arbitrary object $L\in\mathcal{C}_s$ and consider an epimorphism $\alpha:Q^n\lrt L,$
where $Q$ is projective $\Lambda$-module. By part $(ii)$, there exists a $\mathcal{C}_j$-preenvelope
$f:Q^n\lrt K_j$. Since for any $s\geq j$, $\mathcal{C}_s\subseteq\mathcal{C}_j$, $L\in\mathcal{C}_j$
and in particular, there is the following commutative diagram;
{\footnotesize{\[\xymatrix{Q^n\ar[dr]_{\alpha} \ar[rr]^{f^n} &&\ar[dl]^{\psi} K_j^n \\ & L& }\]}}
because $f$ is $\mathcal{C}_j$-preenvelope. We set $Y:=K_j^n$.\\
{\sc step 2}: We show that for any $j$, there exits an object $X_j$ of $\A_j$ such that
there is no any epimorphism $\oplus_{i>j}M_i^{(s_i)}\lrt X_j$
in which $M_i^,$s are objects of $\A_i^,$s and $s_i$ is a set for any $i$. Assume that this is not the case.
By our construction, $\A_j=\{X_1, \cdots, X_t\}$, up to isomorphism. As $\oplus_{i=1}^t X_i$ is finitely generated,
we could assume that there exists a $\Lambda$-epimorphism $\oplus_{i=j+1}^nM_{i}^{m_i}\lrt \oplus_{j=1}^t X_j$
for some positive integers $n, m_i$.
Therefore, Main property (a) of \ref{s1} gives rise to the inequality
${\text {min}}\{\mu^*(X_{j})\mid 1\leq j\leq t\}\geq {\text {min}}\{\mu^*(M_{i})\mid j+1\leq i\leq n\}.$
Since by construction of $\A_i^,$s, none of modules $X_j$ is not a direct summand of $\oplus_{i=j+1}^nM_i$,
the equality may not be accomplished.
Letting $\mu^*(M_{s})={\text {min}}\{\mu^*(M_{i})\mid j+1\leq i\leq n\}$,
by step 1, there is a $\Lambda$-epimorphism
$\oplus_{j=1}^t X_j^{m_j}\lrt M_{s}$, for some integers $m_j>0$, implying that
$\mu^*(M_{s})>\mu^*(X_j)$, for some $1\leq j\leq t$, and so we derive a contradiction.\\
{\sc step 3}:
{As we have seen in step 2, for any $i$, there exists an object $M_i\in\A_i$ such that there does not exist any epimorphism
$\varphi:\oplus_{j>i}M_j^{(s_j)}\lrt M_i$, where $M_j^,$s are objects of $\A_j^,$s and $s_j$ is a set for any $j$.}
Now, for any $i$, we take only one of such modules $M_i$ and denote the class consisting of
all these modules by $\F$. Since $\A_i^,$s are pairwise disjoint infinite family,
$\F$ will be an infinite set of indecomposable pairwise
non-isomorphic $\CM$ $\Lambda$-modules. Therefore, similar to the argument given in
the proof of Theorem \ref{the2}, we get $\Lambda$-homomorphisms $f_i:M_i\lrt S$ such that
for any $j>i$, each composition map $M_j\lrt M_i\st{f_i}\lrt S$ is zero, where $S=\oplus_{j=1}^nS_j$
and $\{S_1, \cdots, S_n\}$ is the complete list of non-isomorphic simple $\Lambda$-modules.
As $\alpha':G'\lrt S$ is a right $\CM$-approximation, one may obtain a $\Lambda$-homomorphism
$g_i:M_i\lrt G'$, for any $i$.
Setting $g=(g_i)_{i\in I}:\oplus_{i\in I}M_i\lrt G'$, {where $I$ is a subset of $\N$}, we have
an exact sequence of $\Lambda$-modules;
$0\lrt K\st{\theta}\lrt \oplus_{i\in I}M_i\oplus \Lambda^n\st{[g~~~\beta]}\lrt G'\lrt 0$.
Clearly, $K$ is $\omega$-Gorenstein projective and so, by the hypothesis, it can be written
as a direct sum of indecomposable
finitely generated modules, say $K=\oplus_{j\in J}K_j$. Now the remainder of the proof
goes along the same lines of the method given in the proof of Theorem \ref{th1},
by replacing $\Lambda$ and $S$ with $R$ and $k$, respectively.
So we omit it.
\end{proof}

Since over Gorenstein algebras, $\omega$-Gorenstein projective $\Lambda$-modules
are just Gorenstein projective modules, as a direct consequence of Theorem \ref{th3}
together with Corollary \ref{cor10}, we recover Chen's theorem \cite[Main Theorem]{che}.

\begin{cor}Let $\Lambda$ be a Gorenstein artin algebra. Then the following conditions are equivalent:
\begin{enumerate}
\item $\Lambda$ is of finite
$\CM$-type.

\item Any Gorenstein projective $\Lambda$-module is fully
decomposable.
\item
 The category of all indecomposable finitely generated Gorenstein projective $\Lambda$-modules is  of bounded length.
 \end{enumerate}\end{cor}

{\bf Acknowledgments.} {The authors are grateful to the referee for reading
the paper very carefully and giving a lot of valuable suggestions kindly and patiently,
especially for pointing out some unclear points in the proof of the previous version of Theorem \ref{fcmt}. 
We also would like to thank Claus Michael Ringel for some helpful conversations.}


\begin{thebibliography}{9999}

\bibitem{af}{\sc F. W. Anderson and K. R. Fuller,} {\sl  Modules with decompositions that complement
direct summands,} J. Algebra {\bf 22} (1972), 241-253.


\bibitem{a1}{\sc M. Auslander,} {\sl A functorial approach to representation theory, in Representations of Algebra,}
Workshop Notes of the Third Inter. Confer., Lecture Notes Math. 944, 105-179, Springer-Verlag, 1982.

\bibitem{a2}{\sc M. Auslander,} {\sl Large modules over Artin algebras,} in: Algebra, Topology,
and category theory (a collection of papers in honor of Samuel Eilenberg), pp. 1-17. Academic Press, New York, 1976.

\bibitem{a11}{\sc M. Auslander,} {\sl Representation theory of Artin algebras. II,}
Comm. Algebra {\bf 1} (1974), 269-310.


\bibitem{ab}{\sc M. Auslander and R. O. Buchweitz,} {\sl The homological theory of
maximal Cohen-Macaulay approximations,} M$\acute{e}$m. Soc. Math. France (N.S.) (1989),
no. 38, 5-37, Colloque en l'honneur de Pierre Samuel (Orsay, 1987).

\bibitem{AR1} {\sc  M. Auslander and I. Reiten,} {\sl Cohen-Macaulay and Gorenstein artin algebras,} in {Representation
theory of finite groups and finite-dimensional algebras (Bielefeld 1991)}, Progress in mathematics, 95 (eds G. O. Michler and C. M. Ringel)
(Birkhauser, Basel, 1991), pp.221-245.

\bibitem{AR} {\sc M. Auslander and I. Reiten,} {\sl Applications of contravariantly finite subcategories,}
Adv. Math. {\bf 86} (1991), no. 1, 111-152.


\bibitem{as1}{\sc M. Auslander and S. O. Smal$\varnothing$,} {\sl Preprojective modules
over artin algebras,} J. Algebra {\bf 66} (1980), 61-122.


\bibitem{be}{\sc A. Beligiannis,} {\sl On algebras of finite Cohen-Macaulay type},
Adv. Math. {\bf 226} (2011), no. 2, 1973-2019.


{\bibitem{bh}{\sc W. Bruns and J. Herzog,} {\sl Cohen-Macaulay Rings,} Cambridge Stud. Adv. Math., Vol. 39,
Cambridge Univ. Press, 1993.}

\bibitem{c}{\sc S. U. Chase,} {\sl Direct product of modules,} Trans. Amer. Math. Soc.
{\bf 97} (1960), 457-473.

\bibitem{ch}{B. Chen,} {\sl The Gabriel-Roiter measure for representation-finite hereidtary algebras},
 J.Algebra. {\bf 309}(2007), 292-317.

\bibitem{che}{\sc X. W. Chen,} {\sl An Auslander-type result for Gorenstein
projective modules,} Adv. Math. {\bf 218} (2008), 2043-2050.


\bibitem{di} {\sc E. Dieterich,} {\sl Representation types of group rings
over complete discrete valuation rings}, Integral representations and applications
(Oberwolfach, 1980,) Lecture Notes in Math., vol. 882, Springer, Berlin, 1981, pp. 369-389.

{\bibitem{ej}{\sc E. E. Enochs and O. M. G. Jenda,} {\sl Relative homological algebra,} de Gruyter Exp. Math., vol. 30, de Gruyter, Berlin, 2000.}


\bibitem{fgr} {\sc R. M. Fossum, P. A. Griffith and I. Reiten,} {\sl Trivial extensions of abelian categories,}
Lecture Notes in Math. VOl. 456. Berlin: Springer-Verlag, 1975.

\bibitem{ga}{\sc P. Gabriel,} {\sl Indecomposable Representations II,}
Symposia Mathematica. Vol. XI, Academic Press, London, 1973, pp. 81-104.

\bibitem{gt}{\sc R. G\"{o}bel and J. Trlifaj,} {\sl Approximations and Endomorphism Algebras of Modules,}
de Gruyter Expositions in Mathematics 41, xxiv + 640 str., W. de Gruyter,
Berlin - New York 2006.

\bibitem{hh}{\sc M. Hochster and C. Huneke,} {\sl Absolute integral closures are big Cohen--Macaulay
algebras in charactristic $\p$}, Bull. Amer. Math. Soc. (N. S.) {\bf 24} (1991), 137-143.

\bibitem{hs}{\sc M. Harada and Y. Sai,} {\sl On the categories of indecomposable
modules, $I$,} Osaka J. Math. {\bf 7} (1970), 323-344.

\bibitem{hoc}{\sc M. Hochster,} {\sl Big and small Cohen-Macaulay modules,}
Module theory (Proc. Special Session, Amer. Math. Soc., Univ. Washington,
Seattle, Wash., 1977), Lecture Notes in Math.,
vol. 700 (Springer, Berlin, 1979), 119-142.

\bibitem{ho}{\sc H. Holm,} {\sl The structure of balanced big Cohen-Macaulay modules
over Cohen-Macaulay rings,} Glasg. Math. J. {\bf 59} (2017), 549-561.

\bibitem{hj} {\sc H. Holm and P. J${\o}$rgensen,} {\sl Semi-dualizing modules and related Gorenstein homological
dimension,} J. Pure Appl. Algebra {\bf 205} (2006), 423-445.

\bibitem{hl}
{\sc C. Huneke; G. Leuschke},  {\sl Two theorems about maximal Cohen-Macaulay modules,}
 Math. Ann. {\bf 324} (2000), no. 2, 391-404.

\bibitem{hr}{\sc P. Hilton and D. Rees,} {\sl Natural maps of extension
functors and a theorem of R. G. Swan,} Proc. Camb. Phil. Soc. {\bf 57} (1961), 489-502.


\bibitem{lz}{\sc Z. W. Li and P. Zhang,} {\sl Gorenstein algebras of finite Cohen-Macaulay type,} Adv. Math. {\bf 223}
(2010), 728-734.

\bibitem{le}{\sc H. Lenzing,} {\sl Homological transfer from finitely presented to infinite modules, Abelian group theory} (Honolulu, Hawaii, 1983), Lecture Notes in Math., vol. 1006, Springer, Berlin, 1983, pp. 734-761.

\bibitem{lw2} {\sc G. J. Leuschke and R. Wiegand, }{\sl Brauer-Thrall Theory for Maximal Cohen-Macaulay Modules,}
 2013, Commutative Algebra, pp. 577-592.

\bibitem{lw1} {\sc G. J. Leuschke and R. Wiegand, }{\sl Cohen-Macaulay Representations,}
Mathematical Surveys and Monographs, American Mathematical Society,
Providence, RI, 2012.

\bibitem{j} {\sc J. P. Jans,} {\sl On the indecomposable representations of algebras,}
Ann. of Math. (2) {\bf 66} (1957), 418-429.

\bibitem{ma} {\sc H. Matsumura,} {\sl Commutative ring theory,} second edition,
Cambridge Studies in Advanced Mathematics, 8, Cambridge University Press, Cambridge, 1989.

\bibitem{mo} {\sc A. Moore,} {\sl The Auslander and Ringel-Tachikawa theorem for
submodule embeddings,} Comm. Algebra {\bf 38} (2010), 3805-3820.

\bibitem{ri3}{\sc C. M. Ringel,} {\sl Foundation of the representation theory of artin algebras,
Using the Gabriel-Roiter measure,} in: Trends in Representation Theory of Algebras and Related Topics,
Workshop Queretaro, Mexico, 2004, in: Contemp. Math., vol. 406, Amer. Math. Soc., 2006, pp. 105-135.




\bibitem{ri2}{\sc M. C. Ringel}, {\sl The Auslander bijections: How morphisms are determined by modules},
Bull. Math. Sci. {\bf 3} (2013), no. 3, 409-484.

\bibitem{ri44}{\sc M. C. Ringel,} {The first Brauer-Thrall conjecture,} Models, modules and abelian groups,
369-374, Walter de Gruyter, Berlin, 2008.

\bibitem{ri1} {\sc C. M. Ringel,} {\sl The Gabriel-Roiter measure,}
 Bull. Sci. Math. {\bf 129} (2005), 726-748.

\bibitem{rt}{\sc C. M. Ringel and H. Tachikawa,} {\sl QF-3 rings,} J. Reine Angew. Math.
{\bf 272} (1975), 49-72.

\bibitem{ro} {\sc A. V. Roiter,} {\sl Unboundedness of the dimension of the indecomposable
representations of an algebra which has infinitely many indecomposable representations,} Izv. Akad. Nauk SSSR.
Ser. Mat. {\bf 32} (1968), 1275-1282.

\bibitem{ru}{\sc W. Rump,} {\sl The category of lattices over a lattice-finite ring,}
 Algebras and Repr. Theory {\bf 8}(2005), 323-345.

\bibitem{ru1}{\sc W. Rump,} {\sl Lattice-finite rings,} Algebras and Repr.
Theory {\bf 8} (2005), 375-395.

\bibitem{sha}{\sc R. Y. Sharp}, {\sl Cohen-Macaulay properties for balanced big
Cohen-Macaulay modules,} In mathematical proceeding of the Cambridge Philosophical
Society, vol. 90, Cambridge Univ. Press, 1981, pp. 229-238.

\bibitem{t}{\sc H. Tachikawa}, {\sl QF-3 rings and categories of projective modules,}
J. Algebra {\bf 28} (1974), 408-413.


\bibitem{yo1}{\sc Y. Yoshino,} {\sl Brauer-Thrall type theorem for maximal Cohen-Macaulay modules,}
J. Math. Soc. Japan {\bf 39} (1987), no. 4, 719-739.

\bibitem{yo}{\sc Y. Yoshino,} {\sl Cohen-Macaulay Modules over Cohen-Macaulay Rings,}
Lond. Math. Soc. Lecture Notes Series, 146, Cambridge Univ. Press, Cambridge, 1990.

\end{thebibliography}
\end{document}